\documentclass[12pt]{article}

\usepackage{amsmath}
\usepackage{amssymb}
\usepackage{amsthm}
\usepackage{latexsym}
\usepackage{cite}
\usepackage{psfrag}
\usepackage{epsfig}
\usepackage{graphicx}
\usepackage[latin1]{inputenc}

\newtheorem{theorem}{Theorem}[section]

\newtheorem{conjecture}[theorem]{Conjecture}
\newtheorem{observation}[theorem]{Observation}
\theoremstyle{definition}

\newtheorem{remark}[theorem]{Remark}

\makeatletter \@addtoreset{equation}{section} \makeatother

\begin{document}
\title{Tables, bounds and graphics of the smallest
known sizes of complete caps in the spaces \\$\mathrm{PG}(3,q)$ and $\mathrm{PG}(4,q)$
\thanks{The research of  D. Bartoli, S. Marcugini and F.
Pambianco was
 supported in part by Ministry for Education, University
and Research of Italy (MIUR) (Project ``Geometrie di Galois e
strutture di incidenza'')
 and by the Italian National Group for Algebraic and Geometric Structures and their Applications
 (GNSAGA - INDAM).
 The research of A.A.~Davydov and A.A. Kreshchuk was carried out at the IITP RAS at the expense of the Russian
Foundation for Sciences (project 14-50-00150).}
}
\date{}
\author{Daniele Bartoli \\
{\footnotesize Dipartimento di Matematica e Informatica,
Universit\`{a}
degli Studi di Perugia, }\\
{\footnotesize Via Vanvitelli~1, Perugia, 06123, Italy.
E-mail: daniele.bartoli@unipg.it}\\
 \and Alexander A.
Davydov,  Alexey A.
Kreshchuk \\
{\footnotesize Institute for Information Transmission Problems
(Kharkevich
institute)}\\ {\footnotesize Russian Academy of Sciences,
 Bol'shoi Karetnyi per. 19, GSP-4, Moscow,
127994}\\{\footnotesize Russian Federation. E-mail: \{adav,krsch\}@iitp.ru}
\and Stefano Marcugini and Fernanda Pambianco \\
{\footnotesize Dipartimento di Matematica e Informatica,
Universit\`{a}
degli Studi di Perugia, }\\
{\footnotesize Via Vanvitelli~1, Perugia, 06123, Italy. E-mail:
\{stefano.marcugini,fernanda.pambianco\}@unipg.it}}

\maketitle

\begin{abstract}
In this paper we present and analyze computational results concerning
small complete caps in the projective spaces $\mathrm{PG}(N,q)$ of dimension $N=3$ and $N=4$
over the finite field of order $q$. The results have been
obtained using randomized greedy
algorithms and the algorithm with fixed order of points (FOP).
The computations have been done in relatively wide regions of $q$ values;
such wide regions are not considered in literature for $N=3,4$.
The new complete caps are the smallest
known.  Basing on them, we obtained new upper bounds on $t_2(N,q)$,
the minimum size of a complete cap
in $\mathrm{PG}(N,q)$, in particular,
\begin{align*}
&t_{2}(N,q)<\sqrt{N+2}\cdot q^{\frac{N-1}{2^{\vphantom{H}}}}\sqrt{\ln q},\quad q\in L_{N},\quad N=3,4,
\displaybreak[2]\\
&t_{2}(N,q)<\left(\sqrt{N+1}+\frac{1.3}{\ln (2q)}\right)q^{\frac{N-1}{2^{\vphantom{H}}}}\sqrt{\ln q},\quad q\in L_{N},\quad N=3,4,
\end{align*}
where
\begin{align*}
&L_{3}:=\{q\le 4673, ~q\  \textrm{prime}\}  \cup \{5003,6007,7001,8009\},\\
&L_{4}:=\{q\le 1361, ~q\  \textrm{prime}\}  \cup \{1409\}.
\end{align*}
Our investigations and results allow to conjecture that these bounds hold for all $q$.
\end{abstract}


\section{Introduction}
Let $\mathrm{PG}(N,q)$ be the $N$-dimensional projective space over the
Galois field $\mathbb{F}_q$ of order $q$. A cap $\mathcal{K}$  in  $\mathrm{PG}(N,q)$ is a
set of points no three of which are collinear.  A cap
$\mathcal{K}$ is complete if it is not contained in a larger
cap or, equivalently, if every point of $\mathrm{PG}(N,q) \setminus
\mathcal{K}$ is collinear with two points of $\mathcal{K}$.
Caps in $\mathrm{PG}(2,q)$ are also called arcs and they have been
widely studied by many authors in the past decades. In
particular, we refer to the surveys and the results
in the works
\cite{BDFKMP-PIT2014,BDFKMP-JG2015-2,BDFKMP-ArXiv2015Greedy,BDFKMP-arXivFOP2015,%
BDFMP-DM,BDFMP-JG2013,BDFMP-JG2015,BFMP-JG2013,DFMP_JG2005,DFMP2009,HirsSt-old,HirsStor-2001,HirsThas-2015,KV,Oster}
(see also the references therein) for the known constructions
and bounds on the size of complete arcs in projective planes.
If $N>2$ only few constructions and bounds are known.

Caps and in particular arcs have been intensively studied for
their connection with Coding Theory. A linear $q$-ary code with length $n$,
dimension $k$, and minimum distance~$d$ is denoted by  $[n,k,d]_{q}$. If a parity-check matrix of a linear $q$-ary code
is obtained by taking as
columns the homogeneous coordinates of the points of a cap in
$\mathrm{PG}(N,q)$, then the code has  minimum distance $4$ (with the exceptions of
the complete 5-cap in $\mathrm{PG}(3,2)$ giving rise to the
$[5,1,5]_{2}$ code and the complete 11-cap in $\mathrm{PG}(4,3)$
corresponding to the Golay $[11,6,5]_{3}$ code). In
particular, complete caps of size $n$ in
$\mathrm{PG}(N,q)$ correspond to non-extendable
 $[n, n-N-1,4]_{q}$ codes. In the case $N=2$ these codes are MDS, that
is they attain the Singleton bound, whereas if $N=3$ they are
Almost MDS, since their Singleton defect is equal to $1$.

Another important parameter concerning linear codes is the
 covering radius. The
covering radius of an $[n,k,d]_{q}$
 code $\mathcal{C}$ is the
minimum integer $r=r(\mathcal{C})$ such that any vector of
$\mathbb{F}^n_{q}$ has distance at most $r$
from $\mathcal{C}$. Complete caps correspond to quasi-perfect
linear codes, that is codes with
$r(\mathcal{C})=\left\lfloor \frac{d-1}{2} \right \rfloor +1$,
since they have minimum distance $4$ and
covering radius $2$; see also
\cite{BLP1998,BPW1989,CHLL1997,CKMS1985,CLS1986}.
The covering density $\mu(\mathcal{C})$, introduced in
\cite{CLS1986}, is one of the parameters characterizing the
covering quality of an $[n, k, d]$-code $\mathcal{C}$ and it is
defined by
$$\mu(\mathcal{C}) =\frac{1}{q^{n-k}}\sum_{i=0}^{r(\mathcal{C})}(q-1)^i\binom{n}{i}.$$

Note also that caps are connected with quantum
codes; see e.g. \cite{CapsQuantCod,Tonchev}.

In general, a central problem concerning caps is to determine
the spectrum of the possible sizes of complete caps in a given
space; see \cite{HirsSt-old,HirsStor-2001} and the references
therein. Of particular interest for applications to Coding
Theory is the lower part of the spectrum; in fact, small
complete caps in projective Galois spaces correspond to
quasi-perfect linear codes with minimum distance 4 and small
density; see for example \cite{CKMS1985,DFMP2009}.

Let $t_2(N,q)$ be \emph{the
minimum size of a complete cap in $\mathrm{PG}(N,q)$}.

The exact values of $t_2(N,q)$ are
known only for small $q$. For instance, $t_2(3,q)$ is known
only for $q\leq 7$; see \cite[Tab. 3]{DFMP2009}.

Whereas the trivial lower bound for $t_2(N,q)$ is $\sqrt 2 q^{(N-1)/2}$,
general constructions of complete caps whose size is close to
this lower bound are only known for $q$ even; see
\cite{DFMP2009,GDT1991,DGMP2010,Giulietti2007,PS1996}.
According to the survey paper \cite{HirsStor-2001}, the
smallest known complete caps in $\mathrm{PG}(3,q)$, with $q$ arbitrary
large, have size approximately $q^{3/2}/2$ and were presented
by Pellegrino in 1998 \cite{Pellegrino1998}. However,
Pellegrino's completeness proof appears to present a major gap,
and counterexamples can be found; see \cite[Sect.
2]{BFG2013}. Recently, using a modification of the approach of
\cite{KV}, the probabilistic upper bound $cq^\frac{N-1}{2}\log^{300} q$, with
$c$ constant,  for the value $t_2(N,q)$ has been obtained; see~\cite{BMP2014Svetlog,BMP2015}.
Computer assisted results
on small complete caps in $\mathrm{PG}(N,q)$ and $\mathrm{AG}(N,q)$ are given in
\cite{BFG2013,DMP_JG2004,DFMP2009,FaPasSch_2012,Oster,Platoni}.
Here and further, $\mathrm{AG}(N,q)$ is the $N$-dimensional affine space over the
field $\mathbb{F}_q$.

In this paper we obtain by computer searches
 results concerning \emph{upper bounds}
on the functions $t_2(3,q)$ and $t_2(4,q)$. These searches requested a huge amount of memory and execution time. In particular, we
constructed small complete caps in $\mathrm{PG}(3,q)$ and $\mathrm{PG}(4,q)$ using two
different approaches\footnote{In this work,
calculations were performed using computational resources of
Multipurpose Computing Complex of National Research Centre
``Kurchatov Institute'', http://computing.kiae.ru}: the algorithm with fixed
order of points (FOP), for $q\in
L_{3}$ in $\mathrm{PG}(3,q)$ and $q\in L_{4}$ in $\mathrm{PG}(4,q)$, and randomized greedy algorithms, for $q\in G_{3}$ in $\mathrm{PG}(3,q)$ and $q\in
G_{4}$ in $\mathrm{PG}(4,q)$, where
\begin{align}
&L_{3}:=\{q\le 4673, ~q\  \textrm{prime}\}  \cup \{5003,6007,7001,8009\},\label{eq1_L3}\\
&G_{3}:=\{q\le 3701,~ q\  \textrm{prime}\} \cup \{3803,3907,4001,4289\},\label{eq1_G3}\\
&L_{4}:=\{q\le 1361, ~q\  \textrm{prime}\}  \cup \{1409\},\label{eq1_L4}\\
&G_{4}:=\{q\le 463, ~q\  \textrm{prime}\}.\label{eq1_G4}
\end{align}

Note that  such relatively wide regions of $q$ values are not considered in
literature for $\mathrm{PG}(3,q)$ and $\mathrm{PG}(4,q)$.

Using the data obtained by
the computer searches we present different functions which
approximate the values $t_2(3,q)$ and
$t_2(4,q)$, as done in
\cite{BDFKMP-PIT2014,BDFMP-DM,BDFKMP-ArXiv2015Greedy,%
BDFKMP-arXivFOP2015,BDFMP-JG2013,BDFKMP-JG2015-2,BDFMP-JG2015,BFMP-JG2013}
for the minimum size of complete arcs in projective planes
$\mathrm{PG}(2,q)$. The main estimates obtained in this paper
 are given in the following theorem, see also Sections \ref{PG3q} and \ref{PG4q}.

\begin{theorem}\label{th:main}
Let $t_2(N,q)$ be the
minimum size of a complete cap in the projective space $\mathrm{PG}(N,q)$.
Let $L_3$ and $L_4$ be the sets of values of $q$ given by relations \eqref{eq1_L3} and \eqref{eq1_L4}, respectively.
The following upper bounds on
$t_2(N,q)$  hold.

\textbf{A.} Upper bounds with the constant multiplier $\sqrt{N+2}$:
\begin{align}
t_{2}(N,q)<\sqrt{N+2}\cdot q^{\frac{N-1}{2^{\vphantom{H}}}}\sqrt{\ln q},\quad q\in L_{N},\quad N=3,4.\label{eq1_th1_const_param}
\end{align}

\textbf{B.} Upper bounds with a decreasing multiplier $\beta_N(q)$:
\begin{align}
t_{2}(N,q)<\beta_N(q)q^{\frac{N-1}{2^{\vphantom{H}}}}\sqrt{\ln q},\quad\beta_N(q)=\sqrt{N+1}+\frac{1.3}{\ln (2q)},\quad q\in L_{N},\quad N=3,4. \label{eq1_bnd_decreas}
\end{align}
\end{theorem}

Our investigations and results (see figures and observations in Sections \ref{PG3q} and \ref{PG4q})
 allow to conjecture that the estimates of Theorem \ref{th:main}, especially the
bound with constant multiplier $\sqrt{N+2}$, hold for every prime power~$q$.
\begin{conjecture}\label{conject}
In $\mathrm{PG}(3,q)$ and $\mathrm{PG}(4,q)$, the upper bounds
\eqref{eq1_th1_const_param}, \eqref{eq1_bnd_decreas} hold for all $q$.
\end{conjecture}

\begin{remark}
\label{remPG2q}
In the works \cite{BDFKMP-JG2015-2,BDFKMP-ArXiv2015Greedy,BDFKMP-arXivFOP2015}, the sizes
of small complete arcs in $\mathrm{PG}(2,q)$ are given for all power prime $q\le301813$.
In this work, we obtained complete arcs
in $\mathrm{PG}(2,q)$ for $301813<q\le321007$, $q$ power prime.
The results of \cite{BDFKMP-JG2015-2,BDFKMP-ArXiv2015Greedy,BDFKMP-arXivFOP2015}
and of this work give the following upper bounds for $\mathrm{PG}(2,q)$:
\begin{align*}
t_{2}(2,q)<1.05\sqrt{3q\ln q}<\sqrt{2+2}\cdot q^\frac{2-1}{2}\sqrt{\ln q}, \quad q\le321007.
\end{align*}
So, the upper bounds \eqref{eq1_th1_const_param} hold also for $N=2$ in a wide region of $q$ values.
\end{remark}

As far as this is known to the authors, complete caps obtained
in this work are the smallest
known in literature  for $\mathrm{PG}(3,q)$ with $q\in\{61,67,71,73,79,83\}$, $97\le q\in L_3$, and
$\mathrm{PG}(4,q)$ with $17\le q\in L_4$. In particular, the results of this work improve ones
of the papers \cite{DFMP2009,FaPasSch_2012,Platoni}.

The paper is organized as follows. In Section \ref{Algorithms},
we describe the main features of the algorithms used in our
searches. In Section \ref{Types of Bounds}, some
{types  of upper bounds} on $t_2(N,q)$ are
discussed. In Section \ref{sec_content}, we shortly give the content of tables collecting
sizes of small complete caps obtained with the help of the algorithms of Section
\ref{Algorithms} (the tables are placed in Appendix). In Sections \ref{PG3q} and \ref{PG4q}, we  analyze
the results presented in the tables and illustrated the analysis by graphics. In Section \ref{sect_conclus}, we do some conclusions from
the present work.

Some results of this paper were briefly
presented in \cite{BDFMP-OC2013}.

\section{Algorithms for small caps in $\mathrm{PG}(N,q)$}
\label{Algorithms}

In this section we describe two different
algorithms used to construct small complete caps in $\mathrm{PG}(3,q)$
and $\mathrm{PG}(4,q)$. First of all note that the number of points of
$\mathrm{PG}(N,q)$ is of order $q^N$ and for instance, if $q\simeq 5\cdot
10^3$ then $|\mathrm{PG}(3,q)|\simeq 1.2 \cdot 10^{11}$: this represents
a strong constraint for any algorithm which investigates
subsets of points in  projective spaces.

\subsection{Algorithm with fixed order of points (FOP)}\label{Subsection:FOP}
This algorithm is a particular type of random algorithm. Some
variants of the algorithm FOP for $\mathrm{PG}(2,q)$ and $\mathrm{PG}(3,q)$ are
given in
\cite{BDFKMP-PIT2014,BDFKMP-JG2015-2,BDFKMP-arXivFOP2015,BDFMP-OC2013,BDFMP-JG2015}.
In this work we describe the algorithm FOP for the arbitrary
space $\mathrm{PG}(N,q)$.

 Firstly we fix a particular order on the points
of $\mathrm{PG}(N,q)$. The algorithm builds a complete cap step by step
adding a new point at each step, until  a complete cap is
obtained.

Let $K^{(i-1)}$ be the cap obtained at the $(i-1)$-th step.
Among the points not lying on bisecants of $K^{(i-1)}$, the
first point in the fixed order is added to $K^{(i-1)}$ to obtain $K^{(i)}$.

Suppose
that the points of $\mathrm{PG}(N,q)$ are ordered as
$A_{1},A_{2},\ldots,A_{\frac{q^{N+1}-1}{q-1}}$.  Consider the empty
set as root of the search and let $K^{(j)}$ be the partial
solution obtained in the $j$-th step, as extension of the root.
We put
\begin{align}
&K^{(0)}=\emptyset, \, K^{(1)}=\{A_{1}\},\,K^{(2)}=\{A_{1},A_{2}\},\, m(1)=2,~K^{(j+1)}=K^{(j)}\cup\{A_{m(j)}\},
\displaybreak[0]\notag\\
& m(j)=\min \Big\{ i \in \left[m(j-1)+1,\frac{q^{N+1}-1}{q-1}\right]\: \Big| \: \nexists \;P,Q \in K^{(j)}  : \, A_{i},P,Q \textrm{ are collinear}\Big\},\notag
\end{align}
i.e. $m(j)$ is the minimum subscript $i$ such that the
corresponding point $A_{i}$ does not lie on a bisecant of $K^{(j)}$. The
process ends when a complete cap is obtained, that is no other points can be added.

We decided to choose a particular order on the points of
$\mathrm{PG}(N,q)$. For seek of simplicity, we considered only $q$
prime.  Let the elements of the field
$\mathbb{F}_{q}=\{0,1,\ldots ,q-1\}$ be treated as integers
modulo $q$. Let the points $A_{i}$ of $\mathrm{PG}(N,q)$ be represented
in homogenous coordinates so that
\begin{align*}
A_{i}=(x_{0}^{(i)},x_{1}^{(i)},\ldots,x_{N}^{(i)}),~ x_{j}^{(i)}\in
\mathbb{F}_{q},
\end{align*}
where the leftmost non-zero element is 1.
The points of $\mathrm{PG}(N,q)$ are sorted according to the lexicographic order
on the $(N+1)$-tuples of their coordinates.
 This order is called a
\emph{lexicographical order of points}. We call \emph{lexicap}
a cap obtained by the algorithm FOP with the lexicographical
order of points. We denote by $t_{2}^{L}(N,q)$ the size of a
complete lexicap in $\mathrm{PG}(N,q)$.  It is important that for
 such a
lexicographical order for prime $q$, \emph{the size $t_{2}^{L}(N,q)$
of a complete lexicap and its set of points depend on $N$ and
$q$ only.}

From a geometrical point of view the lexicographical order of
points is a random order. Clearly, different orders on the
points of $\mathrm{PG}(N,q)$ can determine different size of the
complete cap obtained by the algorithm. Due to our experiences in
similar types of search (see
\cite{BDFKMP-PIT2014,BDFMP-DM,BDFMP-JG2013,BDFKMP-JG2015-2,BDFMP-JG2015,BFMP-JG2013})
we can conjecture that the choice of the order determines only
a small perturbation on the size of the complete caps obtained. For instance,
in \cite[Fig. 11]{BDFMP-JG2015}, sizes of complete arcs in $\mathrm{PG}(2,q)$ obtained
by the algorithm FOP with the lexicographical and the so-called Singer orders of points are compared.
The percentage difference  between the sizes  is approximately
in the interval [-4\%, +4\%] for $q\ge1000$.

Connections of the algorithm FOP with algorithms of Coding Theory are noted in
\cite[Remark 3.1]{BDFKMP-JG2015-2} and \cite[Remark 2.1]{BDFMP-JG2015}.

\subsection{Randomized greedy algorithms}\label{Subsection:RandomizedGreedyAlgorithms}

A different approach can be used to obtain small complete caps
in $\mathrm{PG}(N,q)$. In general, small complete caps in $\mathrm{PG}(N,q)$
obtained using randomized gre\-edy algorithms have size smaller
than those obtained with programs of type FOP as described in
the previous subsection; see
\cite{BDFKMP-PIT2014,BDFKMP-ArXiv2015Greedy,BDFKMP-arXivFOP2015,BDFMP-DM,%
BDFMP-JG2013,BDFKMP-JG2015-2,BDFMP-JG2015,BFMP-JG2013,DFMP_JG2005,DMP_JG2004}.

The main difference  between the two types of algorithm is
that at every step a randomized greedy  algorithm maximizes an
objective function $f$ and only some steps are executed in a
random manner. The number of these steps, their ordinal numbers{,}
and some other parameters of the algorithm have been taken
intuitively. Also, if the same maximum of $f$ can be obtained
in distinct ways, one way is chosen randomly.

We start constructing a complete  cap by using a starting point
set $ S_{0} $. In the $i$-th step one point is added to the set
$S_{i-1}$ and we obtain a point set $S_{i}$. As the value of
the objective function $f$ we consider the number of covered
points in $\mathrm{PG}(N,q),$ that is, points that lie on bisecants of
$S_{i}$.

On every  random $i$-th step we take $ d_{q,i}$ randomly chosen
points of\linebreak $\mathrm{PG}(N,q)$ not covered by $ S_{i-1}$ and
compute the objective function $f$ adding each of these $
d_{q,i} $ points to $S_{i-1}$. The point providing the maximum
of $f$ is included into~$S_{i}.$ On every  non-random $j$-th
step we consider all points not covered by $S_{j-1}$ and add to
$S_{j-1}$ the point providing the maximum of $f.$

As $S_{0}$ we can use a subset of points of an arc obtained in previous stages of the search.

A generator of random numbers is used for  random choices.  To
obtain caps with distinct sizes, starting conditions of the
generator are changed for the same set $S_{0}$. In this way the
algorithm works in a convenient limited region of the search
space to obtain examples improving the size of the cap from
which the fixed points have been taken.

In order to obtain arcs with new sizes, sufficiently  many
attempts should be made with  randomized greedy algorithms.
``Predicted'' sizes  could be useful for understanding if a
good result has been obtained.  If the result is not close to
the predicted size, the attempts are continued.

We obtain small complete caps in $\mathrm{PG}(N,q)$ in two
stages.

At the 1-st stage, we take the frame as $S_{0}$ and create  a
starting complete cap $K_{0}$ using in the beginning of the
process $\delta _{q}$  random steps with distinct $d_{q,i}.$
All the subsequent steps are non-random.

At the 2-nd stage we execute $n_{q}$ attempts to get a complete
cap. For every attempt, the starting conditions of the random
generator are different from the previous ones, whereas the set
$S_{0}$ is the same. Two or three among the first five steps of
every attempt are  random, the rest of them are  non-random.

The values  $d_{q,i}$, $\delta _{q}$, and $n_{q}$  are given
intuitively depending on $q$ and (for $d_{q,i})$ on $|S_{i-1}|$
and on the stage of the process. Of course, CPU performance
affects the algorithm parameters choice.

Cap sizes obtained by the randomized greedy algorithms depend
on  many factors, but in general the results are better than
the ones obtained by the algorithm FOP. Unfortunately, this
approach requires a huge amount of execution time and therefore
this type of search has been executed only for a relatively
small region of values of $q$.

\section{General types of bounds for the value $t_2(N,q)$}\label{Types of Bounds}
Let $t_2(N,q)$ be the size of the smallest complete cap in $\mathrm{PG}(N,q)$.
In this section we propose different types of bounds for these values, generalizing the approach proposed for estimates on
 $t_{2}(2,q)$ in \cite{BDFMP-JG2013,BDFKMP-JG2015-2,BDFMP-JG2015}. Also, let
  $t_{2}^{G}(N,q)$  denote the smallest size of a complete cap
  in $\mathrm{PG}(N,q)$ obtained using greedy  algorithms. Finally, remind that
$t_{2}^{L}(N,q)$ is  the size of the complete lexicap in $\mathrm{PG}(N,q)$
obtained by the algorithm FOP with the lexicographical
order of points.

Let   $\beta_{N}(q)$, $\beta_{N}^{G}(q)$, and
$\beta_{N}^{L}(q)$ be some functions of $q$ defined as follows:
\begin{align}
\beta_{N}(q)=\frac{t_{2}(N,q)}{q^{\frac{N-1}{2}}\sqrt{\ln q}},\quad
\beta_{N}^{G}(q)=\frac{t_{2}^{G}(N,q)}{q^{\frac{N-1}{2}}\sqrt{\ln q}},\quad
\beta_{N}^{L}(q)=\frac{t_{2}^{L}(N,q)}{q^{\frac{N-1}{2}}\sqrt{\ln q}}.\label{eq_bNfrac}
\end{align}

From  \eqref{eq_bNfrac} we obtain
\begin{align}
&t_{2}(N,q)=\beta_{N}(q)q^{\frac{N-1}{2}}\sqrt{\ln q}, \quad t_{2}^{G}(N,q)=\beta_{N}^{G}(q)q^{\frac{N-1}{2}}\sqrt{\ln q},\label{eq_bN}\displaybreak[0] \\
&t_{2}^{L}(N,q)=\beta_{N}^{L}(q)q^{\frac{N-1}{2}}\sqrt{\ln q}.\notag
\end{align}

Clearly $t_{2}(N,q)\le \min\{t_{2}^{G}(N,q), t_{2}^{L}(N,q)\}$
and in general,  due to the main features of the two
algorithms, $t_{2}^{G}(N,q)\leq  t_{2}^{L}(N,q)$ always holds.
This implies
\begin{align}
 \beta_{N}(q)\le \beta_{N}^{G}(q)\le \beta_{N}^{L}(q),\quad \beta_{N}(q)\le\min\{\beta_{N}^{G}(q),\beta_{N}^{L}(q)\}.\label{eq_fdb<}
\end{align}

We consider \textbf{two types of
upper bounds} on $t_{2}(N,q)$.

A. \emph{Upper bounds with the constant multiplier $\sqrt{N+2}$}. For
this type,
 we consider upper bounds on   $\beta_{N}(q)$ equal
 to a value dependent on $N$ but independent of $q$.

B. \emph{Upper bounds with a decreasing multiplier $\beta_N(q)$}.
For this type, we find upper bounds on  $\beta_{N}(q)$ as
a decreasing function of $q$ denoted by   $\beta_{N}^{up}(q)$.
This function looks like
\begin{align*}
\beta_{N}^{up}(q)=a+\frac{b}{\ln(cq)},
\end{align*}
where  $a$ is a value dependent on $N$ but independent of $q$, whereas $b,c$ are constants independent of $N$ and $q$.

\section{The content of tables}\label{sec_content}
Results of our computer searches are collected in tables given in Appendix.

In Table 1, for $q\in L_{3}$,
we collected the sizes $t_2^L(3,q)$ ($t_2^L$ for short)
of complete lexicaps in $\mathrm{PG}(3,q)$ obtained using the algorithm FOP with the lexicographical order of points, see
Section~\ref{Subsection:FOP}.

In
Table 2, for $q\in G_{3}$, the sizes
$t_{2}^{G}(3,q)$ ($t_{2}^{G}$ for short) of complete caps in $\mathrm{PG}(3,q)$ obtained using randomized
greedy algorithms, see Section
\ref{Subsection:RandomizedGreedyAlgorithms}, are given.

Note that for
$q\in\{61,67,71,73,79,83,97\}$ sizes $t_{2}^{G}(3,q)$ in Table 2
 improve the ones from \cite[Table 7]{DFMP2009}. Also, the values of $t_{2}^{G}(3,q)$ given in Table 2
are smaller than the sizes of complete caps in $\mathrm{AG}(3,q)$ obtained in \cite[Section 3]{Platoni}.
The improvements are written in Table 2 in bold font.

In Table 3 we collected the sizes
$t_{2}^{L}(4,q)$ ($t_{2}^{L}$ for short) of complete lexicaps
in $\mathrm{PG}(4,q)$, $q\in L_{4}$, obtained by the algorithm
FOP with lexicographical order of points, see Section~\ref{Subsection:FOP}.

In Table 4 we give the sizes
$t_{2}^{G}(4,q)$ ($t_{2}^{G}$ for short) of complete
caps in $\mathrm{PG}(4,q)$, $q\in G_{4}$, obtained by
the randomized greedy algorithms, see Section \ref{Subsection:RandomizedGreedyAlgorithms}.

Note that the size $t_2^G(4,17)$ in Table 4  improves the one from \cite[Table 8]{DFMP2009}.
Also, the values of $t_{2}^{G}(4,q)$ given in Table 4
are smaller than the sizes of complete caps in $\mathrm{AG}(4,q)$ obtained in \cite[Theorem 1.1]{FaPasSch_2012}.
The improvements are written in Table 4 in bold font.

\section{Small complete caps in $\mathrm{PG}(3,q)$}\label{PG3q}

\begin{figure}[htbp]
\includegraphics[width=\textwidth]{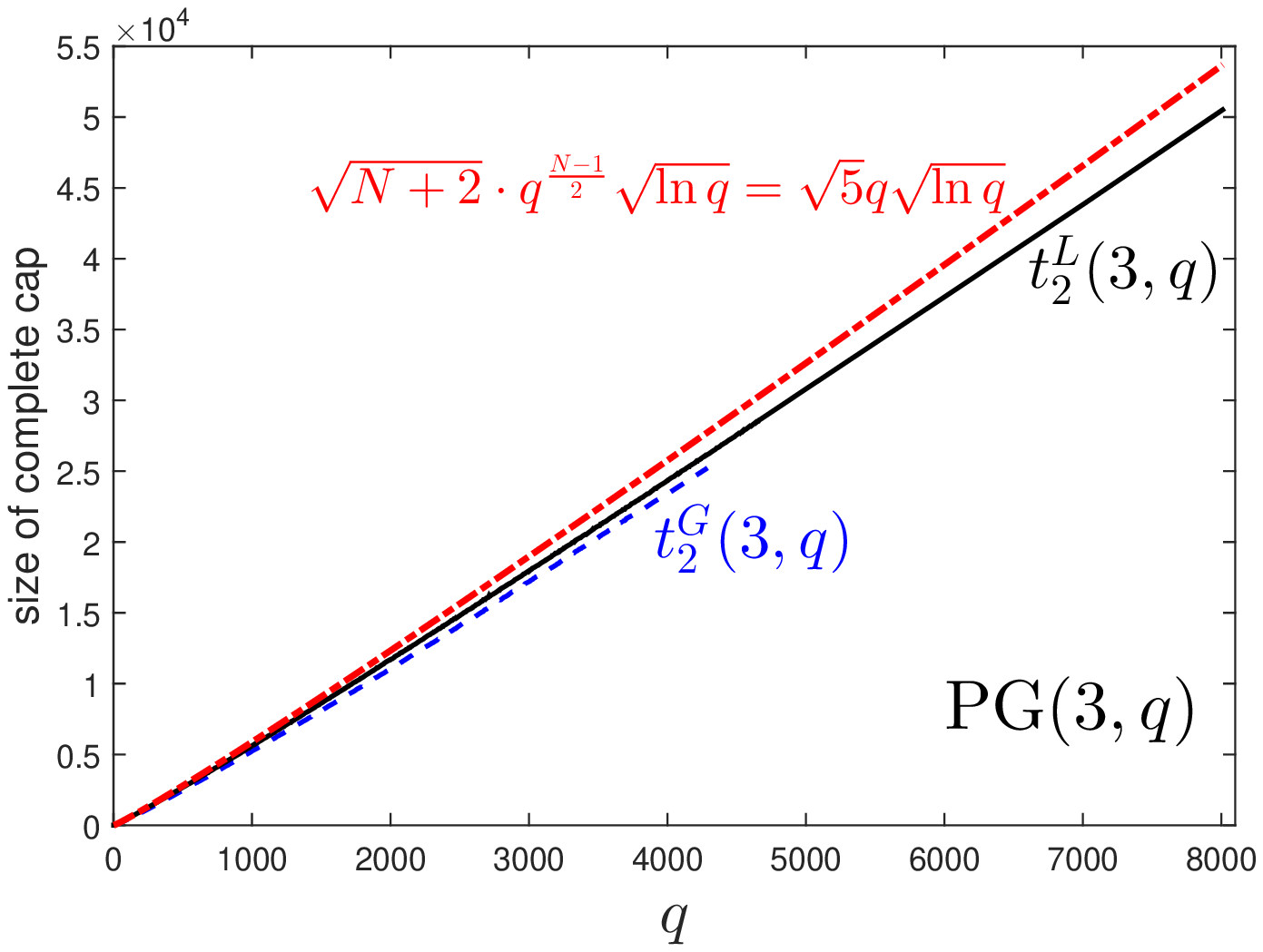}
\caption{$\mathrm{PG}(3,q)$. \textbf{Upper bound} $t_{2}(3,q)
<\sqrt{3+2}\cdot q^{\frac{3-1}{2}}\sqrt{\ln q}=\sqrt{5}q\sqrt{\ln q}$  (\emph{top dashed-dotted red curve})
 \textbf{vs sizes}  $t_{2}^{L}(3,q)$ of complete lexicaps, $q\in L_{3}$ (\emph{the 2-nd solid black curve}) and \textbf{sizes}
$t_{2}^{G}(3,q)$ of complete caps obtained by greedy algorithms, $q\in G_{3}$ (\emph{bottom dashed blue curve}).}
\label{fig_3q_1}
\end{figure}

\begin{figure}[htbp]
\includegraphics[width=\textwidth]{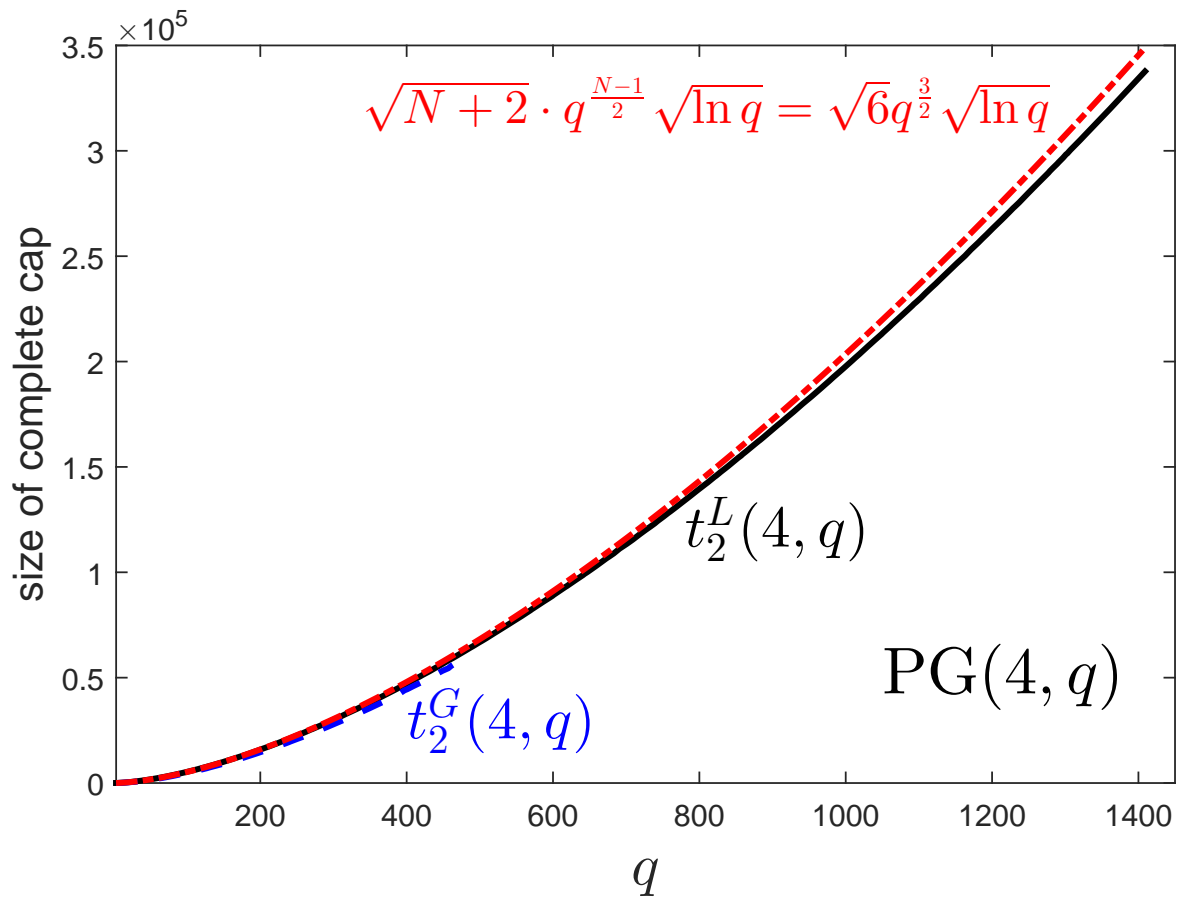}
\caption{$\mathrm{PG}(4,q)$. \textbf{Upper bound}  $t_{2}(4,q)
<\sqrt{4+2}\cdot q^{\frac{4-1}{2}}\sqrt{\ln q}=\sqrt{6}q^\frac{3}{2}\sqrt{\ln q}$  (\emph{top dashed-dotted red curve})
 \textbf{vs sizes}  $t_{2}^{L}(4,q)$ of complete lexicaps, $q\in L_{4}$ (\emph{the 2-nd solid black curve}) and \textbf{sizes}
$t_{2}^{G}(4,q)$ of complete caps obtained by greedy algorithms, $q\in G_{4}$ (\emph{bottom dashed blue curve}).}
\label{fig_4q_1}
\end{figure}

The values $t_2^L(3,q)$ written in Table 1 are shown in
Figure \ref{fig_3q_1} by the 2-nd solid black curve. In turn,  the values $t_2^L(4,q)$, given in Table 3,
 are shown by the 2-nd solid black curve in Figure~\ref{fig_4q_1}.

 The values $t_2^G(3,q)$ from Table 2 are shown in Figure \ref{fig_3q_1}
by the bottom dashed blue curve. In the scale of Figure \ref{fig_3q_1} the curves $t_{2}^{L}(3,q)$
and $t_{2}^{G}(3,q)$ are very close to each other.

Note that for all $q\in G_{3}$, we have
$t_2^G(3,q)< t_2^L(3,q)$. So, as already pointed out  above,
the use of greedy algorithms provides better results, that is
the size of the complete caps obtained is smaller. However,
randomized greedy algorithms require in general more execution
time than algorithm FOP, since they require more investigations
at each step, trying to maximize a particular objective
function as illustrated in Section
\ref{Subsection:RandomizedGreedyAlgorithms}. For this reason we
have been able to obtain the data for a smaller
region of values of $q$ than by the FOP algorithm.

\begin{figure}[htbp]
\includegraphics[width=\textwidth]{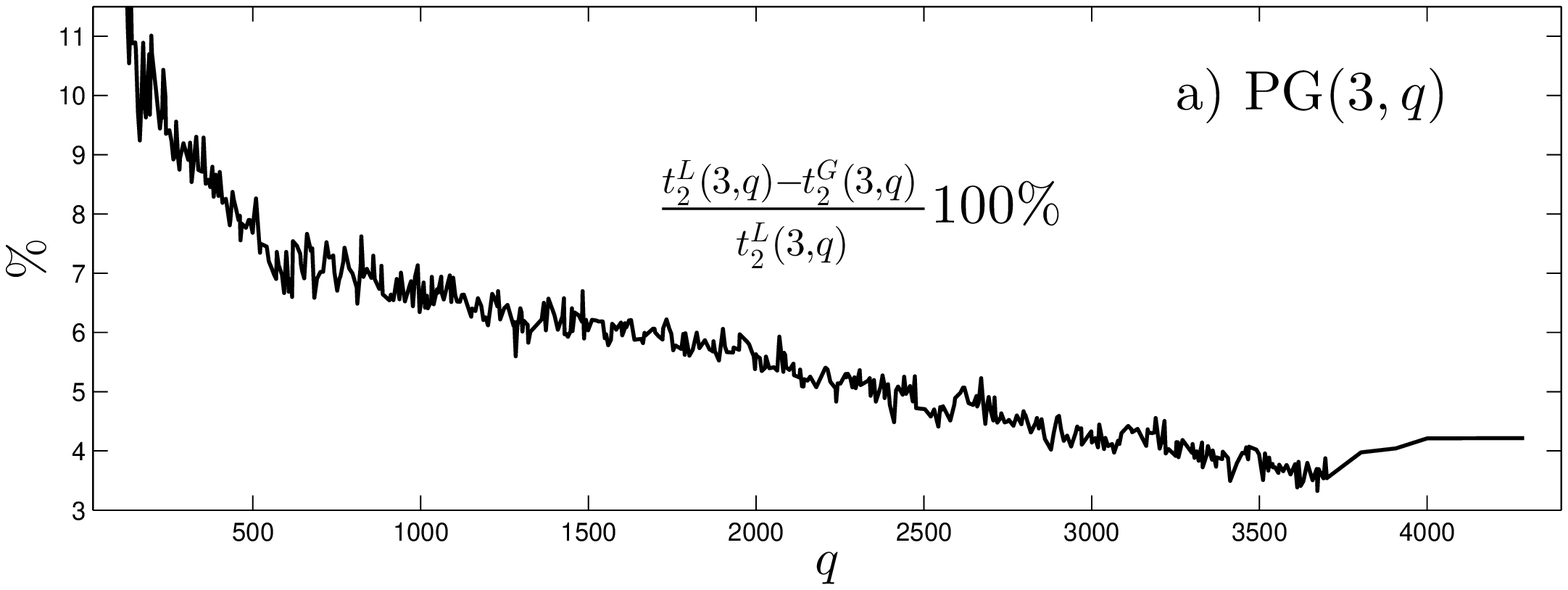}
\includegraphics[width=\textwidth]{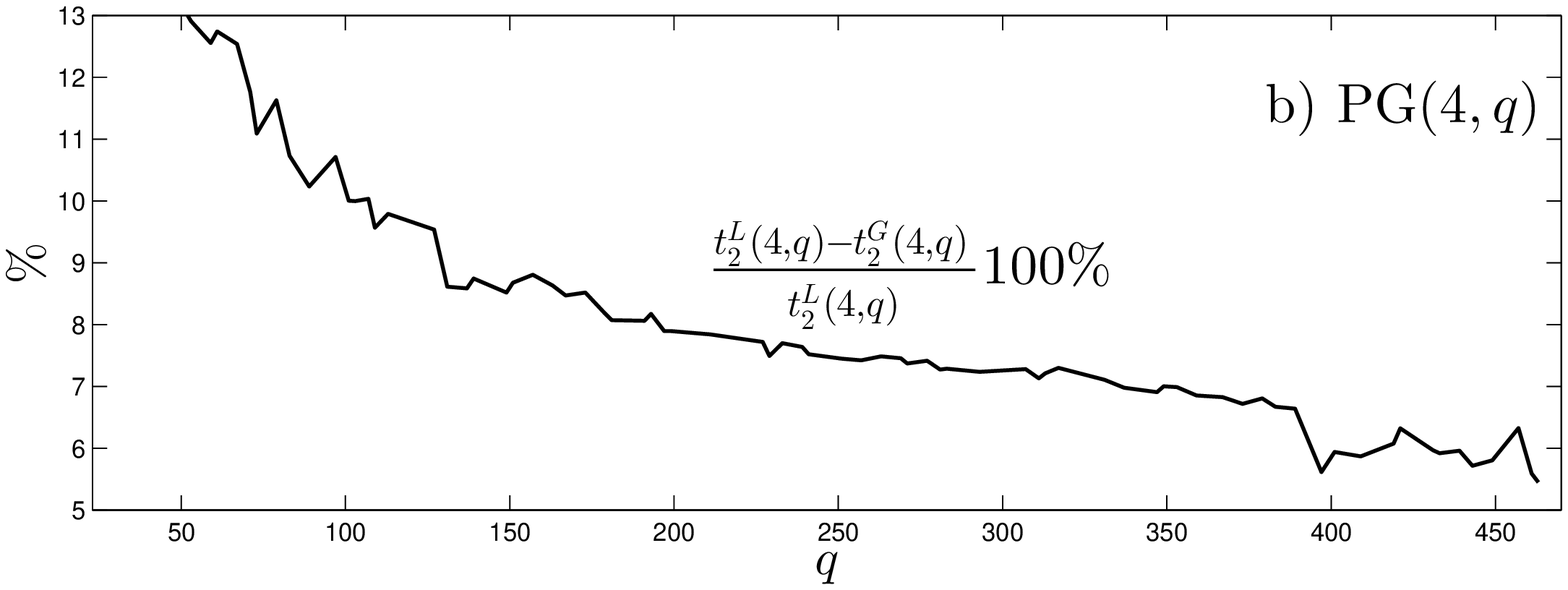}
\caption{\textbf{Percentage difference between  $t_2^L(N,q)$ and
$t_2^G(N,q)$.}
a) $N=3$, $\mathrm{PG}(3,q)$; b) $N=4$, $\mathrm{PG}(4,q)$}
\label{fig_3q4q_2}
\end{figure}

Figure \ref{fig_3q4q_2}a shows the percentage difference between  $t_2^L(3,q)$ and
 $t_2^G(3,q)$.

\begin{observation}\label{remPG3q}
From Tables \emph{1} and \emph{2} and Figure \emph{\ref{fig_3q4q_2}a}, one sees that the percentage difference between  $t_2^L(3,q)$ and
 $t_2^G(3,q)$ given~by
$$\frac{t_{2_{\vphantom{H}}}^{L}(3,q)-t_{2}^{G}(3,q)}{t_{2}^{L^{\vphantom{H}}}(3,q)}100\%$$
is relatively small and it tends to decrease when $q$ grows. In particular, in the region $q\in[503\ldots3701]$ this difference decreases
approximately from 7\% to 4\%.
\end{observation}

\begin{figure}[htbp]
\includegraphics[width=\textwidth]{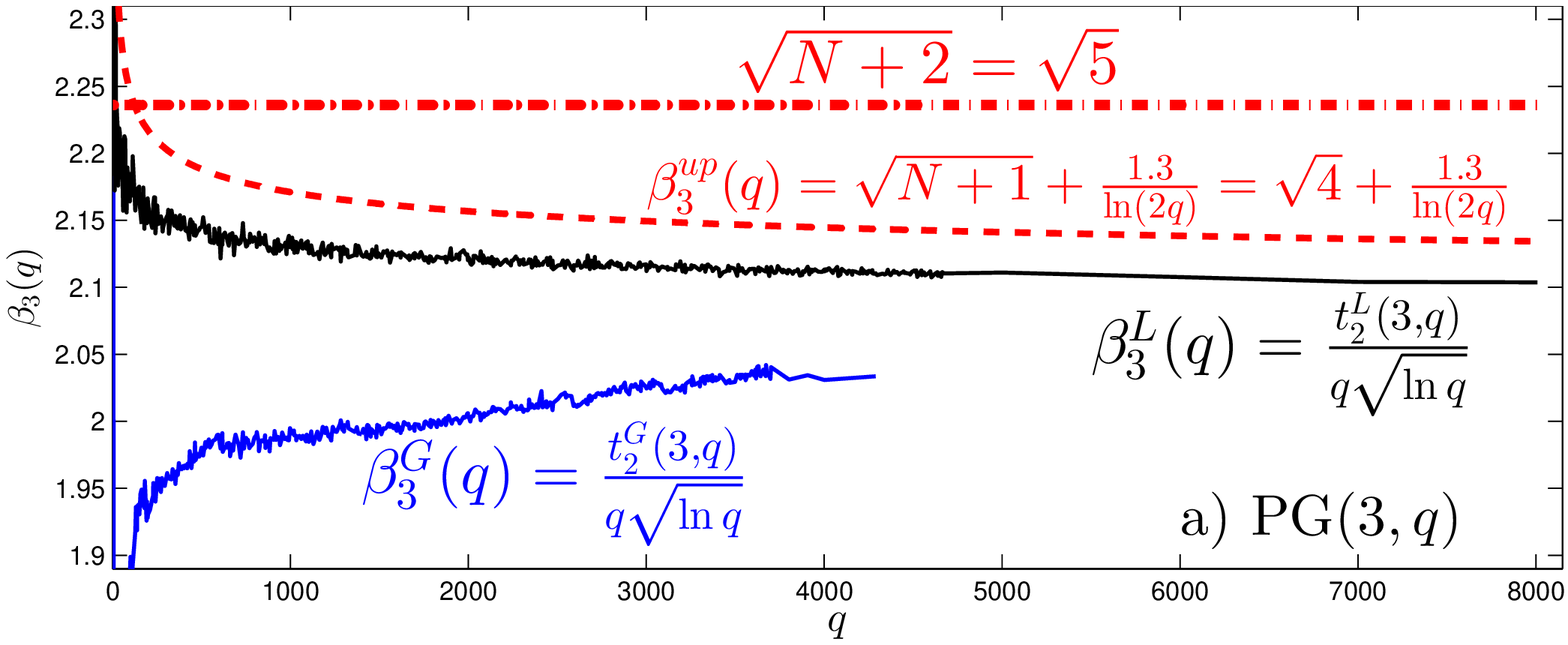}
\includegraphics[width=\textwidth]{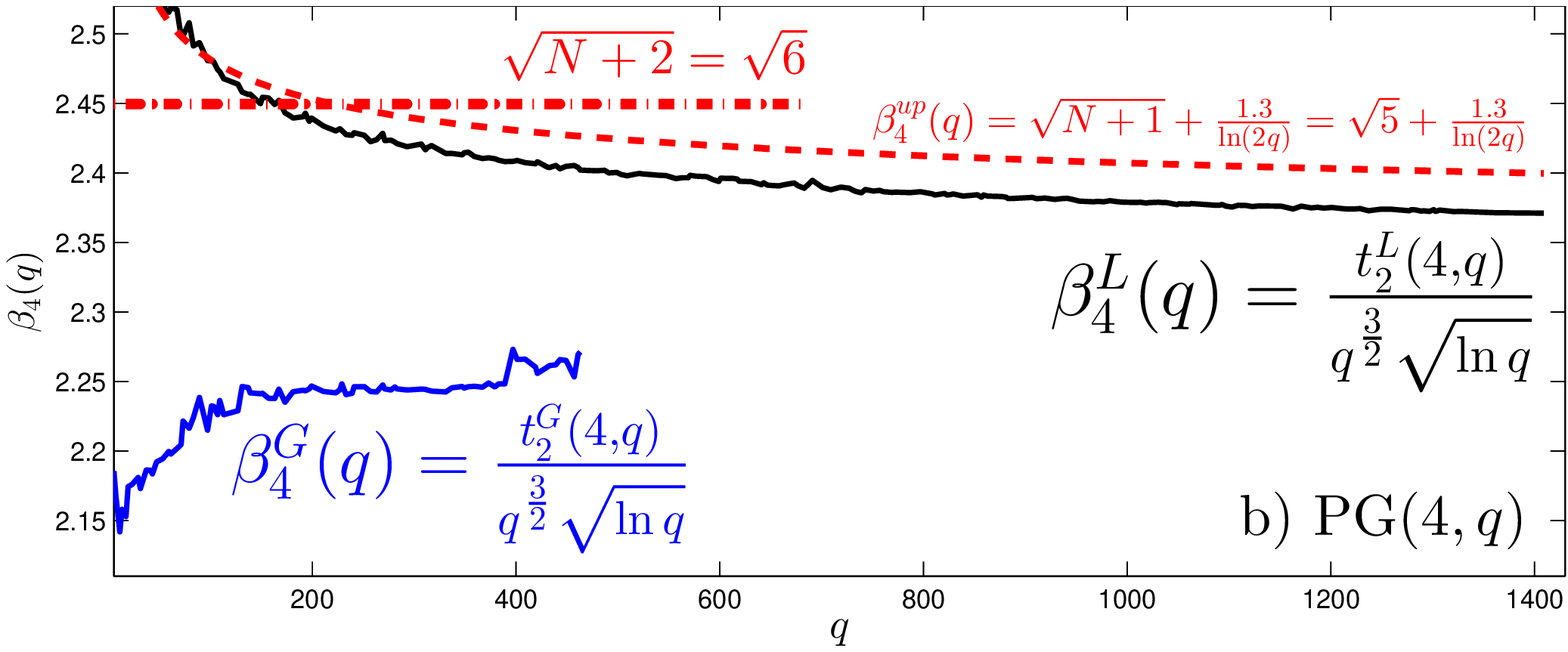}
\caption{\textbf{Upper bounds} $\beta_N(q)=\frac{t_{2}(N,q)}{q^{\frac{N-1}{2}}\sqrt{\ln q}}
<\sqrt{N+2}$
(\emph{dashed-dotted red line $y=\sqrt{N+2}$}) and
$\beta_N(q)<\beta_{N}^{up}(q)=\sqrt{N+1}+\frac{1.3}{\ln (2q)}$
(\emph{top dashed red curve}) \textbf{vs values} of $
\beta_N^L(q)$, $q\in L_{N}$ (\emph{the 2-nd solid
black curve}) and  $\beta_N^G(q)$,
$q\in G_{N}$ (\emph{bottom solid blue curve}). a) $N=3$, $\mathrm{PG}(3,q)$; b) $N=4$, $\mathrm{PG}(4,q)$}
\label{fig_3q4q_3}
\end{figure}

Figure \ref{fig_3q4q_3}a shows the values
 $\beta_{3}^{L}(q)$ and $\beta_{3}^{G}(q)$
 obtained by \eqref{eq_bNfrac} from the sizes collected in Tables
 1 and 2. Also, in this figure, upper bounds
  $\beta^{up}_{3}(q)=\sqrt{N+1}+\frac{1.3}{\ln (2q)}=\sqrt{3+1}+\frac{1.3}{\ln (2q)}=2+\frac{1.3}{\ln (2q)}$ and the line-bound
  $y=\sqrt{N+2}=\sqrt{5}$ are presented in red color.

 By Tables 1,
2 and Figure \ref{fig_3q4q_3}a, it holds that, see \eqref{eq_fdb<},
\begin{align}
&\beta_{3}(q)\le\min\{\beta_{3}^{G}(q),\beta_{3}^{L}(q)\}<\sqrt{N+2}=\sqrt{5},\quad q\in L_{3}\label{eq_PG3q_sqrt5};\displaybreak[3]\\
&\beta_{3}(q)\le\min\{\beta_{3}^{G}(q),\beta_{3}^{L}(q)\}<\beta^{up}_{3}(q)=\sqrt{N+1}+\frac{1.3}{\ln (2q)}=\sqrt{3+1}+\frac{1.3}{\ln (2q)}=\label{eq_PG3q_b3up}\\
&\phantom{\beta_{3}(q)\le\min\{\beta_{3}^{G}(q),\beta_{3}^{L}(q)\}<\beta^{up}_{3}(q)=~}2+\frac{1.3}{\ln (2q)},\quad q\in L_{3}.\notag
\end{align}
This implies upper bounds for $\mathrm{PG}(3,q)$ in Theorem \ref{th:main}.

The upper bound  \eqref{eq1_th1_const_param} for $N=3$, based on \eqref{eq_PG3q_sqrt5},
is shown by the dashed-dotted red curve in Figure \ref{fig_3q_1}.
This bound is
presented also by the dashed-dotted red line $y=\sqrt{N+2}=\sqrt{5}$
in Figure \ref{fig_3q4q_3}a. The bound \eqref{eq1_bnd_decreas} for $N=3$, based on
\eqref{eq_PG3q_b3up}, is given
by the dashed red curve in Figure \ref{fig_3q4q_3}a.

\begin{figure}[htbp]
\includegraphics[width=\textwidth]{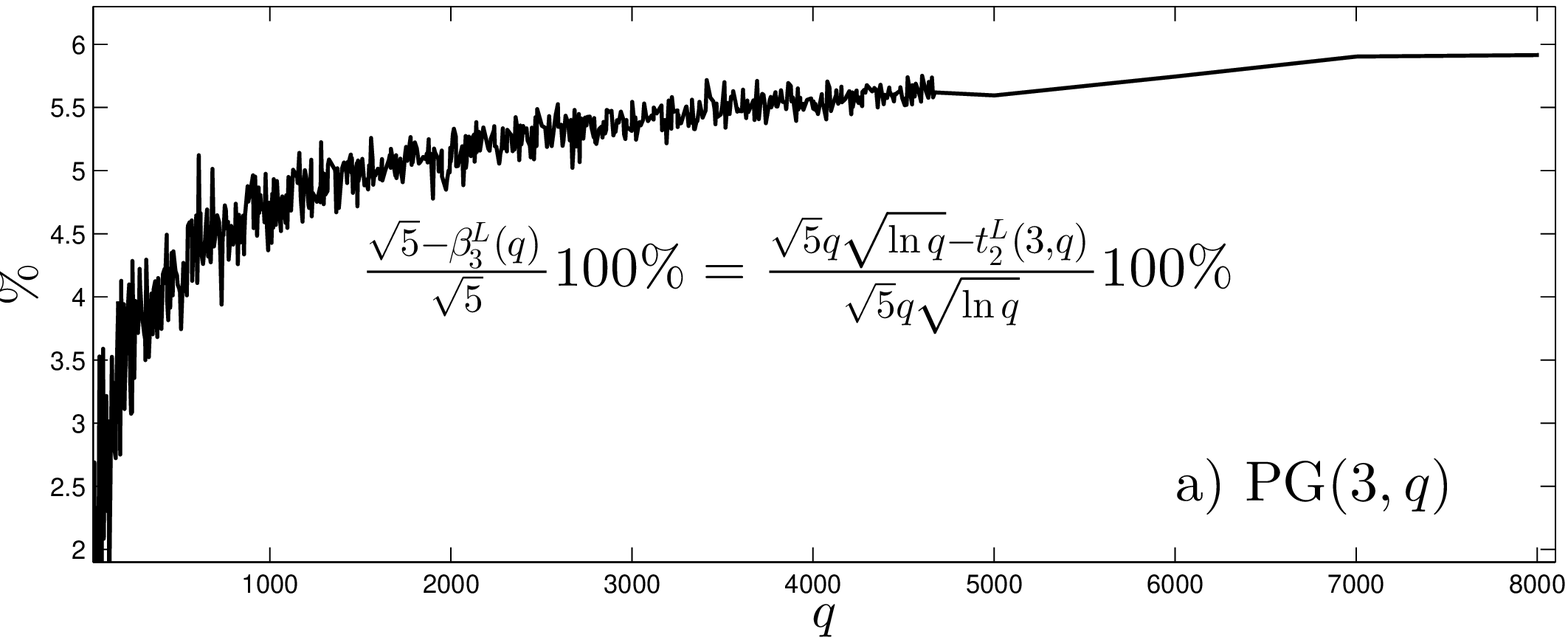}
\includegraphics[width=\textwidth]{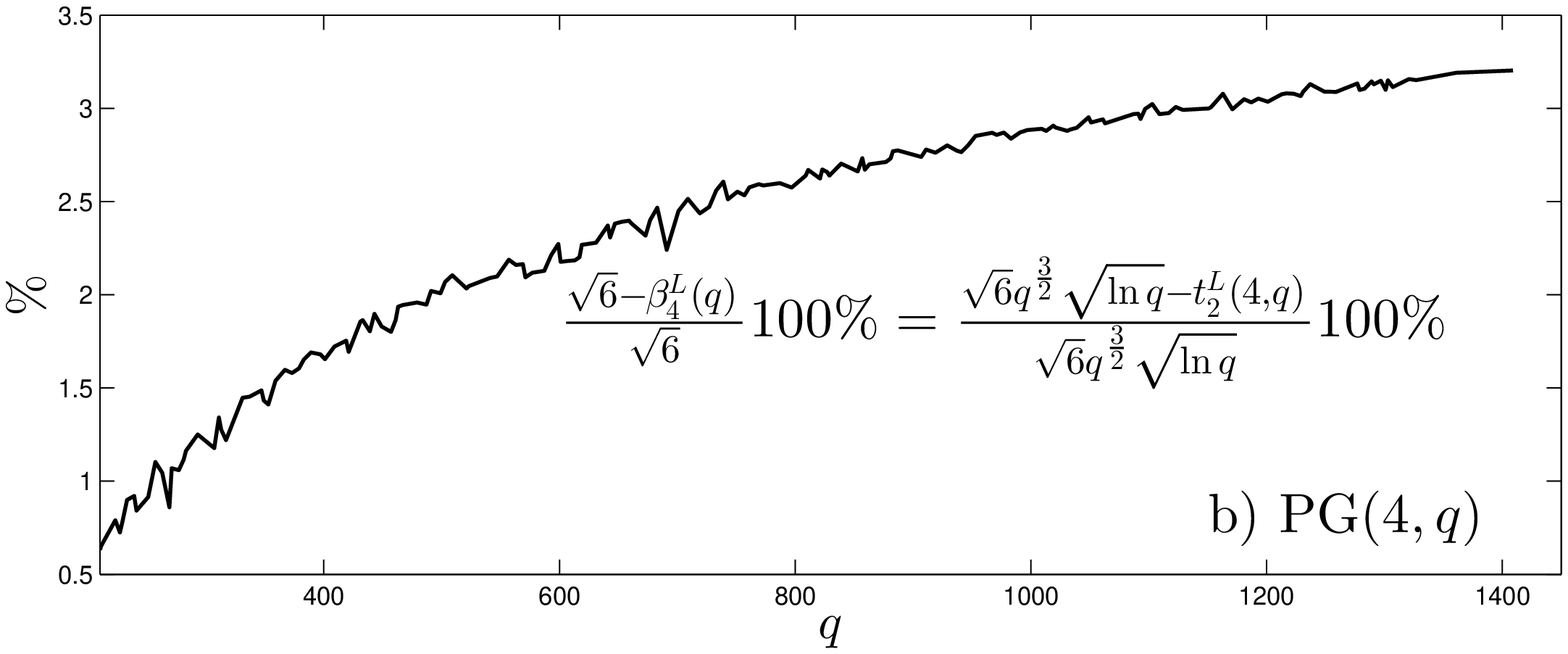}
\caption{
\textbf{
Percentage difference between
$\sqrt{N+2}\cdot
q^{\frac{N-1}{2}}\sqrt{\ln q}$
and $t_2^L(N,q)$.
}
a) $N=3$, $\mathrm{PG}(3,q)$; b) $N=4$, $\mathrm{PG}(4,q)$}
\label{fig_3q4q_4}
\end{figure}

Figure \ref{fig_3q4q_4}a shows the percentage differences between  $\sqrt{5}$ and
 $\beta_{3}^{L}(q)$ and $\sqrt{5}q\sqrt{\ln q}$ and $t_{2}^{L}(3,q)$. (These percentage differences are equal to each other.)

\begin{observation}\label{observ_PG3q}
From Table \emph{1} and Figure \emph{\ref{fig_3q4q_3}a} one sees that the curve
 $\beta_{3}^{L}(q)$ has a decreasing trend.
Therefore the difference
$\sqrt{5}-\beta_{3}^{L}(q)$, and  the corresponding percent
differences
\begin{align*}
\frac{\sqrt{5}-\beta_{3}^{L}(q)}{\sqrt{5}}100\%=\frac{\sqrt{5}q\sqrt{\ln q}-t_{2}^{L}(3,q)}{\sqrt{5}q\sqrt{\ln q}}100\%,
\end{align*}
tend to increase when $q$ grows, see Figure \emph{\ref{fig_3q4q_4}a}. This raises confidence in the correctness of the
bound \eqref{eq1_th1_const_param} for  $N=3$.
\end{observation}

Concerning the execution time, the search for the small
complete cap in $\mathrm{PG}(3,7001)$ lasted 2 months with  a processor
AMD Opteron(TM) Processor 6212, $2.6$ Ghz, and used $85$GB of
memory.

\section{Small complete caps in $\mathrm{\mathrm{PG}}(4,q)$}\label{PG4q}

Figures \ref{fig_4q_1} and \ref{fig_3q4q_3}b show the values
$t_{2}^{L}(4,q)$, $t_{2}^{G}(4,q)$ collected in Tables
3,  4  and the
corresponding values  $\beta_{4}^{L}(q)$, $\beta_{4}^{G}(q)$ obtained by
\eqref{eq_bNfrac}. Also, in these figures, upper bounds
are presented. Note that in the scale of Figure \ref{fig_4q_1} the curves $\sqrt{6}q^\frac{3}{2}\sqrt{\ln q}$,
$t_{2}^{L}(4,q)$, and $t_{2}^{G}(4,q)$ are very closed to each other.

Figure \ref{fig_3q4q_2}b shows the percentage difference between  $t_2^L(4,q)$ and
 $t_2^G(4,q)$.

\begin{observation}\label{observPG4q}
Even if for all $q\in G_{4}$ the inequality
$t_{2}^{G}(4,q)<t_{2}^{L}(4,q)$ holds, see Figure
\emph{\ref{fig_4q_1}}, the difference in percentage between these two
values given by
$$\frac{t_{2_{\vphantom{H}}}^{L}(4,q)-t_{2}^{G}(4,q)}{t_{2}^{L^{\vphantom{H}}}(4,q)}100\%$$
is relatively small and it tends to decrease when $q$ grows, see Figure \emph{\ref{fig_3q4q_2}b}.  In particular, in the region $q\in[101\ldots443]$ this difference decreases
approximately from 10\% to 5\%.
\end{observation}

Figure \ref{fig_3q4q_3}b shows the values
 $\beta_{4}^{L}(q)$ and $\beta_{4}^{G}(q)$
 obtained by \eqref{eq_bNfrac} from the sizes collected in Tables
 3 and 4. Also, in this figure, upper bounds
  $\beta^{up}_{4}(q)=\sqrt{N+1}+\frac{1.3}{\ln (2q)}=\sqrt{4+1}+\frac{1.3}{\ln (2q)}=\sqrt{5}+\frac{1.3}{\ln (2q)}$ and the line-bound
  $y=\sqrt{N+2}=\sqrt{6}$ are presented in red color.

  From Tables 3,
4 and Figure \ref{fig_3q4q_3}b, we have, see \eqref{eq_fdb<},
\begin{align}
&\beta_{4}(q)\le\min\{\beta_{4}^{G}(q),\beta_{4}^{L}(q)\}<\sqrt{N+2}=\sqrt{6},\quad q\in L_{4}; \label{eq_PG4q_sqrt6}\displaybreak[1]\\
&\beta_{4}(q)\le\min\{\beta_{4}^{G}(q),\beta_{4}^{L}(q)\}<\beta^{up}_{4}(q)=
\sqrt{N+1}+\frac{1.3}{\ln (2q)}=\sqrt{4+1}+\frac{1.3}{\ln (2q)}=\label{eq_PG4q_b4}\\
&\phantom{\beta_{4}(q)\le\min\{\beta_{4}^{G}(q),\beta_{4}^{L}(q)\}<\beta^{up}_{4}(q)=~}\sqrt{5}+\frac{1.1}{\ln q},\quad q\in L_{4}. \notag
\end{align}

This implies bounds for $\mathrm{PG}(4,q)$ in Theorem \ref{th:main}.

The upper bounds \eqref{eq1_th1_const_param} for $N=4$, based on
\eqref{eq_PG4q_sqrt6}, are shown
by the dashed-dotted red curves in Figure \ref{fig_4q_1}. This bound
is presented also by the dashed-dotted red line  $y=\sqrt{N+2}=\sqrt{6}$
in Figure \ref{fig_3q4q_3}b. The bound \eqref{eq1_bnd_decreas} for $N=4$, based on
\eqref{eq_PG4q_b4}, is given
by the dashed red curve in Figure \ref{fig_3q4q_3}b.

Figure \ref{fig_3q4q_4}b shows the percentage differences between  $\sqrt{6}$ and
 $\beta_{4}^{L}(q)$ and $\sqrt{6}q\sqrt{\ln q}$ and $t_{2}^{L}(4,q)$. (These percentage differences are equal to each other.)

\begin{observation}\label{observ_PG4q}
From Table \emph{3} and Figure \emph{\ref{fig_3q4q_3}b} one sees that the curve
$\beta_{4}^{L}(q)$ have a decreasing trend.
Therefore the difference
$\sqrt{6}-\beta_{4}^{L}(q)$ and  the corresponding percent
differences
\begin{align*}
\frac{\sqrt{6}-\beta_{4}^{L}(q)}{\sqrt{6}}100\%=\frac{\sqrt{6}q^\frac{3}{2}\sqrt{\ln q}-t_{2}^{L}(3,q)}{\sqrt{6}q^\frac{3}{2}\sqrt{\ln q}}100\%
\end{align*}
tend to increase when $q$ grows, see Figure \emph{\ref{fig_3q4q_4}b}. This raises confidence in the correctness of the
bound \eqref{eq1_th1_const_param} for  $N=4$.
\end{observation}

\section{Conclusion} \label{sect_conclus}
In this paper we presented and analyze computational results
concerning small complete caps in $\mathrm{PG}(3,q)$, $q\le
4673$, $q$ prime, and $q=5003,6007,7001$, $8009$, and
$\mathrm{PG}(4,q)$, $q\leq 1361$, $q$ prime, and $q=1409$.

The results have been obtained using randomized greedy
algorithms and the algorithm with fixed order of points (FOP).
Tables 1--4 and Figures
\ref{fig_3q_1}--\ref{fig_3q4q_4} show that the sizes $t_2^G(N,q)$ of complete
caps obtained by greedy algorithms are smaller than sizes
$t_2^L(N,q)$ of complete
caps
formed by the algorithm FOP with the lexicographical order of points. This allows, in particular, to
increase the regions of $q$ values where the proposed upper
bounds hold, see Figures
\ref{fig_3q_1}, \ref{fig_3q4q_3} and relations \eqref{eq_fdb<}, \eqref{eq_PG3q_sqrt5}, \eqref{eq_PG3q_b3up},
 \eqref{eq_PG4q_sqrt6}, \eqref{eq_PG4q_b4}.

In the other side, the percent difference between
$t_2^L(N,q)$ and $t_2^G(N,q)$ is relatively small and it decreases when $q$ grows, see Observations \ref{remPG3q}
and \ref{observPG4q}. Execution time
of greedy algorithms is essentially greater than for the algorithm FOP.
The sizes $t_2^G(N,q)$ depend not only on $q$ and $N$ but also on
parameters $d_{q,i}$, $\delta _{q}$, $n_{q}$ of greedy algorithms,
see Section \ref{Subsection:RandomizedGreedyAlgorithms}. These parameters
are not always chosen optimal due to restrictions of the computer time.

At the same time,
the sizes $t_2^L(N,q)$ depend on $q$ and $N$ only. Therefore the behavior of the curves
  $\beta_{N}^{L}(q)$ obtained from $t_2^L(N,q)$
  allows to understand the order of value and effectively estimate  the smallest sizes
$t _2(N,q)$ of complete caps for $N=3,4$ in the considered regions of $q$, see Figures \ref{fig_3q4q_3}, \ref{fig_3q4q_4}.

Moreover, the decreasing trend of the curves
 $\beta_{N}^{L}(q)$, see Figures \ref{fig_3q4q_3}, \ref{fig_3q4q_4}, allow us to conjecture that the upper bounds on
$t_2(3,q)$ and $t_2(4,q)$ we obtained, especially  the bounds
\eqref{eq1_th1_const_param} with constant multiplier $\sqrt{N+2}$,
 hold for any $q$ prime power.

As far as this is known to the authors, new complete caps obtained
in this work are the smallest
known in literature.

\section{Appendix. Tables of sizes of the small complete caps in $\mathrm{PG}(3,q)$ and $\mathrm{PG}(4,q)$}
In Table 1, for $q\in L_{3}$,
we collected the sizes $t_2^L(3,q)$ ($t_2^L$ for short)
of complete lexicaps in $\mathrm{PG}(3,q)$ obtained using the algorithm FOP with the lexicographical order of points, see
Section~\ref{Subsection:FOP}.

In Table 2, for $q\in G_{3}$, the sizes
$t_{2}^{G}(3,q)$ ($t_{2}^{G}$ for short) of complete caps in $\mathrm{PG}(3,q)$ obtained using randomized
greedy algorithms, see Section
\ref{Subsection:RandomizedGreedyAlgorithms}, are given.

In Table 3 we collected the sizes
$t_{2}^{L}(4,q)$ ($t_{2}^{L}$ for short) of complete lexicaps
in $\mathrm{PG}(4,q)$, $q\in L_{4}$, obtained by the algorithm
FOP with lexicographical order of points, see Section~\ref{Subsection:FOP}.

In Table 4 we give the sizes
$t_{2}^{G}(4,q)$ ($t_{2}^{G}$ for short) of complete
caps in $\mathrm{PG}(4,q)$, $q\in G_{4}$, obtained by
the randomized greedy algorithms, see Section \ref{Subsection:RandomizedGreedyAlgorithms}.
\newpage
\begin{center}
\textbf{Table 1.} Sizes $t_2^L(3,q)=t_2^L$ of complete lexicaps in
$\mathrm{PG}(3,q)$, {$q\in L_{3}$}\medskip\\
\begin{tabular}{@{}r@{\,\,\,}r@{\,\,}|@{\,\,}r@{\,\,\,}r@{\,\,}|@{\,\,}r@{\,\,\,}r@{\,\,}|@{\,\,}
r@{\,\,\,}r@{\,\,}|@{\,\,}r@{\,\,\,}r@{\,\,}|@{\,\,}r@{\,\,\,}r@{\,\,}|@{\,\,}r@{\,\,\,}r@{}}
\hline
$q$&$t_2^L$&$q$&$t_2^L$&$q$&$t_2^L$&$q$&$t_2^L$&$q$&$t_2^L$&$q$&$t_2^L$&$q$
&$t_{2_{\vphantom{H_H}}}^{L^{\vphantom{H^H}}}$\\
\hline
 2 & 8  &  3 & 8  &  5 & 16  &  7 & 23  &  11 & 37  &  13 & 49  &  17 & 69  \\
 19 & 71  &  23 & 91  &  29 & 118  &  31 & 125  &  37 & 156  &  41 & 175  &  43 & 183  \\
 47 & 202  &  53 & 232  &  59 & 257  &  61 & 273  &  67 & 304  &  71 & 324  &  73 & 328  \\
 79 & 356  &  83 & 382  &  89 & 410  &  97 & 449  &  101 & 474  &  103 & 481  &  107 & 502  \\
 109 & 512  &  113 & 540  &  127 & 603  &  131 & 626  &  137 & 660  &  139 & 671  &  149 & 725  \\
 151 & 732  &  157 & 761  &  163 & 790  &  167 & 814  &  173 & 854  &  179 & 874  &  181 & 893  \\
 191 & 944  &  193 & 951  &  197 & 981  &  199 & 990  &  211 & 1050  &  223 & 1112  &  227 & 1138  \\
 229 & 1154  &  233 & 1179  &  239 & 1212  &  241 & 1208  &  251 & 1275  &  257 & 1300  &  263 & 1334  \\
 269 & 1368  &  271 & 1381  &  277 & 1412  &  281 & 1429  &  283 & 1442  &  293 & 1501  &  307 & 1582  \\
 311 & 1605  &  313 & 1619  &  317 & 1628  &  331 & 1720  &  337 & 1750  &  347 & 1801  &  349 & 1813  \\
 353 & 1841  &  359 & 1868  &  367 & 1912  &  373 & 1952  &  379 & 1989  &  383 & 2002  &  389 & 2043  \\
 397 & 2083  &  401 & 2113  &  409 & 2149  &  419 & 2205  &  421 & 2219  &  431 & 2267  &  433 & 2291  \\
 439 & 2329  &  443 & 2342  &  449 & 2384  &  457 & 2419  &  461 & 2447  &  463 & 2462  &  467 & 2476  \\
 479 & 2549  &  487 & 2596  &  491 & 2621  &  499 & 2668  &  503 & 2692  &  509 & 2735  &  521 & 2791  \\
 523 & 2801  &  541 & 2913  &  547 & 2931  &  557 & 2988  &  563 & 3024  &  569 & 3057  &  571 & 3084  \\
 577 & 3112  &  587 & 3170  &  593 & 3195  &  599 & 3248  &  601 & 3254  &  607 & 3260  &  613 & 3317  \\
 617 & 3334  &  619 & 3356  &  631 & 3430  &  641 & 3482  &  643 & 3493  &  647 & 3512  &  653 & 3543  \\
 659 & 3592  &  661 & 3601  &  673 & 3676  &  677 & 3693  &  683 & 3706  &  691 & 3777  &  701 & 3832  \\
 709 & 3873  &  719 & 3934  &  727 & 3992  &  733 & 4044  &  739 & 4056  &  743 & 4080  &  751 & 4117  \\
 757 & 4154  &  761 & 4184  &  769 & 4229  &  773 & 4266  &  787 & 4337  &  797 & 4403  &  809 & 4468  \\
 811 & 4471  &  821 & 4544  &  823 & 4565  &  827 & 4578  &  829 & 4582  &  839 & 4652  &  853 & 4725  \\
 857 & 4764  &  859 & 4769  &  863 & 4784  &  877 & 4856  &  881 & 4886  &  883 & 4897  &  887 & 4920  \\
 907 & 5030  &  911 & 5076  &  919 & 5102  &  929 & 5187  &  937 & 5214  &  941 & 5249  &  947 & 5277  \\
 953 & 5318  &  967 & 5404  &  971 & 5433  &  977 & 5447  &  983 & 5507  &  991 & 5566  &  997 & 5580  \\
 1009 & 5671  &  1013 & 5672  &  1019 & 5724  &  1021 & 5713  &  1031 & 5779  &  1033 & 5807  &  1039 & 5822  \\
 1049 & 5900  &  1051 & 5908  &  1061 & 5979  &  1063 & 5976  &  1069 & 6010  &  1087 & 6132  &  1091 & 6142  \\
 1093 & 6158  &  1097 & 6195  &  1103 & 6205  &  1109 & 6242  &  1117 & 6307  &  1123 & 6332  &  1129 & 6358  \\
 1151 & 6495  &  1153 & 6510  &  1163 & 6554  &  1171 & 6628  &  1181 & 6701  &  1187 & 6716  &  1193 & 6745  \\
 1201 & 6797  &  1213 & 6886  &  1217 & 6912  &  1223 & 6935  &  1229 & 6957  &  1231 & 7002  &  1237 & 7025  \\
 1249 & 7100  &  1259 & 7164  &  1277 & 7265  &  1279 & 7286  &  1283 & 7274  &  1289 & 7343  &  1291 & 7341  \\
 1297 & 7397  &  1301 & 7411  &  1303 & 7415  &  1307 & 7452  &  1319 & 7523  &  1321 & 7519  &  1327 & 7552  \\
 \hline
\end{tabular}
\end{center}

\newpage
\begin{center}
\textbf{Table 1. Continue 1.} Sizes $t_2^L(3,q)=t_2^L$ of complete lexicaps in
$\mathrm{PG}(3,q)$, {$q\in L_{3}$}\medskip\\
\begin{tabular}{@{}r@{\,\,\,}r@{\,\,}|@{\,\,}r@{\,\,\,}r@{\,\,}|@{\,\,}r@{\,\,\,}r@{\,\,}|@{\,\,}
r@{\,\,\,}r@{\,\,}|@{\,\,}r@{\,\,\,}r@{\,\,}|@{\,\,}r@{\,\,\,}r@{\,\,}|@{\,\,}r@{\,\,\,}r@{}}
\hline
$q$&$t_2^L$&$q$&$t_2^L$&$q$&$t_2^L$&$q$&$t_2^L$&$q$&$t_2^L$&$q$&$t_2^L$&$q$
&$t_{2_{\vphantom{H_H}}}^{L^{\vphantom{H^H}}}$\\
\hline
 1361 & 7766  &  1367 & 7830  &  1373 & 7838  &  1381 & 7913  &  1399 & 7997  &  1409 & 8054  &  1423 & 8147  \\
 1427 & 8178  &  1429 & 8176  &  1433 & 8207  &  1439 & 8236  &  1447 & 8281  &  1451 & 8318  &  1453 & 8320  \\
 1459 & 8381  &  1471 & 8441  &  1481 & 8495  &  1483 & 8538  &  1487 & 8530  &  1489 & 8557  &  1493 & 8579  \\
 1499 & 8613  &  1511 & 8676  &  1523 & 8769  &  1531 & 8814  &  1543 & 8895  &  1549 & 8914  &  1553 & 8947  \\
 1559 & 8955  &  1567 & 9020  &  1571 & 9064  &  1579 & 9103  &  1583 & 9131  &  1597 & 9212  &  1601 & 9240  \\
 1607 & 9276  &  1609 & 9294  &  1613 & 9315  &  1619 & 9340  &  1621 & 9369  &  1627 & 9403  &  1637 & 9460  \\
 1657 & 9575  &  1663 & 9608  &  1667 & 9643  &  1669 & 9640  &  1693 & 9810  &  1697 & 9818  &  1699 & 9849  \\
 1709 & 9902  &  1721 & 9969  &  1723 & 9990  &  1733 & 10063  &  1741 & 10101  &  1747 & 10127  &  1753 & 10157  \\
 1759 & 10206  &  1777 & 10297  &  1783 & 10358  &  1787 & 10368  &  1789 & 10393  &  1801 & 10450  &  1811 & 10517  \\
 1823 & 10605  &  1831 & 10653  &  1847 & 10748  &  1861 & 10837  &  1867 & 10869  &  1871 & 10902  &  1873 & 10906  \\
 1877 & 10918  &  1879 & 10943  &  1889 & 11006  &  1901 & 11122  &  1907 & 11113  &  1913 & 11154  &  1931 & 11272  \\
 1933 & 11296  &  1949 & 11376  &  1951 & 11413  &  1973 & 11563  &  1979 & 11591  &  1987 & 11632  &  1993 & 11671  \\
 1997 & 11672  &  1999 & 11698  &  2003 & 11718  &  2011 & 11746  &  2017 & 11778  &  2027 & 11848  &  2029 & 11890  \\
 2039 & 11952  &  2053 & 12027  &  2063 & 12082  &  2069 & 12159  &  2081 & 12184  &  2083 & 12229  &  2087 & 12247  \\
 2089 & 12250  &  2099 & 12297  &  2111 & 12392  &  2113 & 12387  &  2129 & 12482  &  2131 & 12523  &  2137 & 12536  \\
 2141 & 12544  &  2143 & 12588  &  2153 & 12629  &  2161 & 12699  &  2179 & 12787  &  2203 & 12957  &  2207 & 12985  \\
 2213 & 13025  &  2221 & 13062  &  2237 & 13153  &  2239 & 13161  &  2243 & 13208  &  2251 & 13259  &  2267 & 13378  \\
 2269 & 13388  &  2273 & 13405  &  2281 & 13448  &  2287 & 13464  &  2293 & 13524  &  2297 & 13539  &  2309 & 13628  \\
 2311 & 13616  &  2333 & 13760  &  2339 & 13811  &  2341 & 13829  &  2347 & 13849  &  2351 & 13890  &  2357 & 13890  \\
 2371 & 13997  &  2377 & 14066  &  2381 & 14073  &  2383 & 14063  &  2389 & 14123  &  2393 & 14142  &  2399 & 14164  \\
 2411 & 14247  &  2417 & 14280  &  2423 & 14332  &  2437 & 14427  &  2441 & 14466  &  2447 & 14468  &  2459 & 14582  \\
 2467 & 14581  &  2473 & 14654  &  2477 & 14634  &  2503 & 14818  &  2521 & 14934  &  2531 & 15021  &  2539 & 15045  \\
 2543 & 15058  &  2549 & 15121  &  2551 & 15140  &  2557 & 15184  &  2579 & 15283  &  2591 & 15345  &  2593 & 15387  \\
 2609 & 15485  &  2617 & 15563  &  2621 & 15582  &  2633 & 15625  &  2647 & 15723  &  2657 & 15809  &  2659 & 15805  \\
 2663 & 15820  &  2671 & 15934  &  2677 & 15945  &  2683 & 15933  &  2687 & 15985  &  2689 & 15992  &  2693 & 16045  \\
 2699 & 16062  &  2707 & 16090  &  2711 & 16180  &  2713 & 16139  &  2719 & 16172  &  2729 & 16255  &  2731 & 16275  \\
 2741 & 16310  &  2749 & 16371  &  2753 & 16385  &  2767 & 16473  &  2777 & 16552  &  2789 & 16641  &  2791 & 16627  \\
 2797 & 16699  &  2801 & 16713  &  2803 & 16739  &  2819 & 16820  &  2833 & 16905  &  2837 & 16946  &  2843 & 16962  \\
 2851 & 17029  &  2857 & 17062  &  2861 & 17059  &  2879 & 17176  &  2887 & 17225  &  2897 & 17324  &  2903 & 17357  \\
 2909 & 17395  &  2917 & 17429  &  2927 & 17496  &  2939 & 17556  &  2953 & 17657  &  2957 & 17704  &  2963 & 17749  \\
\hline
\end{tabular}
\end{center}

\newpage
\begin{center}
\textbf{Table 1. Continue 2.} Sizes $t_2^L(3,q)=t_2^L$ of complete lexicaps in
$\mathrm{PG}(3,q)$, {$q\in L_{3}$}\medskip\\
\begin{tabular}{@{}r@{\,\,\,}r@{\,\,}|@{\,\,}r@{\,\,\,}r@{\,\,}|@{\,\,}r@{\,\,\,}r@{\,\,}|@{\,\,}
r@{\,\,\,}r@{\,\,}|@{\,\,}r@{\,\,\,}r@{\,\,}|@{\,\,}r@{\,\,\,}r@{\,\,}|@{\,\,}r@{\,\,\,}r@{}}
\hline
$q$&$t_2^L$&$q$&$t_2^L$&$q$&$t_2^L$&$q$&$t_2^L$&$q$&$t_2^L$&$q$&$t_2^L$&$q$
&$t_{2_{\vphantom{H_H}}}^{L^{\vphantom{H^H}}}$\\
\hline
 2969 & 17759  &  2971 & 17782  &  2999 & 17966  &  3001 & 17935  &  3011 & 18024  &  3019 & 18073  &  3023 & 18132  \\
 3037 & 18200  &  3041 & 18220  &  3049 & 18276  &  3061 & 18344  &  3067 & 18357  &  3079 & 18454  &  3083 & 18480  \\
 3089 & 18524  &  3109 & 18643  &  3119 & 18713  &  3121 & 18707  &  3137 & 18818  &  3163 & 18971  &  3167 & 19036  \\
 3169 & 19034  &  3181 & 19109  &  3187 & 19155  &  3191 & 19210  &  3203 & 19218  &  3209 & 19274  &  3217 & 19353  \\
 3221 & 19327  &  3229 & 19402  &  3251 & 19540  &  3253 & 19580  &  3257 & 19586  &  3259 & 19613  &  3271 & 19689  \\
 3299 & 19841  &  3301 & 19853  &  3307 & 19903  &  3313 & 19962  &  3319 & 19981  &  3323 & 20011  &  3329 & 20059  \\
 3331 & 20033  &  3343 & 20159  &  3347 & 20185  &  3359 & 20246  &  3361 & 20261  &  3371 & 20305  &  3373 & 20305  \\
 3389 & 20432  &  3391 & 20447  &  3407 & 20537  &  3413 & 20523  &  3433 & 20681  &  3449 & 20790  &  3457 & 20855  \\
 3461 & 20907  &  3463 & 20885  &  3467 & 20909  &  3469 & 20946  &  3491 & 21082  &  3499 & 21118  &  3511 & 21153  \\
 3517 & 21272  &  3527 & 21292  &  3529 & 21323  &  3533 & 21345  &  3539 & 21389  &  3541 & 21394  &  3547 & 21409  \\
 3557 & 21476  &  3559 & 21476  &  3571 & 21571  &  3581 & 21654  &  3583 & 21661  &  3593 & 21730  &  3607 & 21820  \\
 3613 & 21819  &  3617 & 21872  &  3623 & 21911  &  3631 & 21987  &  3637 & 21993  &  3643 & 22053  &  3659 & 22143  \\
 3671 & 22211  &  3673 & 22228  &  3677 & 22286  &  3691 & 22303  &  3697 & 22397  &  3701 & 22439  &  3709 & 22463  \\
 3719 & 22507  &  3727 & 22560  &  3733 & 22613  &  3739 & 22691  &  3761 & 22789  &  3767 & 22821  &  3769 & 22854  \\
 3779 & 22909  &  3793 & 22979  &  3797 & 23032  &  3803 & 23071  &  3821 & 23164  &  3823 & 23168  &  3833 & 23255  \\
 3847 & 23346  &  3851 & 23385  &  3853 & 23396  &  3863 & 23411  &  3877 & 23507  &  3881 & 23563  &  3889 & 23607  \\
 3907 & 23769  &  3911 & 23768  &  3917 & 23769  &  3919 & 23824  &  3923 & 23832  &  3929 & 23865  &  3931 & 23897  \\
 3943 & 23967  &  3947 & 24005  &  3967 & 24114  &  3989 & 24253  &  4001 & 24351  &  4003 & 24349  &  4007 & 24391  \\
 4013 & 24419  &  4019 & 24472  &  4021 & 24481  &  4027 & 24530  &  4049 & 24631  &  4051 & 24665  &  4057 & 24705  \\
 4073 & 24838  &  4079 & 24798  &  4091 & 24888  &  4093 & 24919  &  4099 & 24989  &  4111 & 25051  &  4127 & 25136  \\
 4129 & 25139  &  4133 & 25184  &  4139 & 25232  &  4153 & 25311  &  4157 & 25312  &  4159 & 25358  &  4177 & 25468  \\
 4201 & 25606  &  4211 & 25720  &  4217 & 25712  &  4219 & 25756  &  4229 & 25802  &  4231 & 25807  &  4241 & 25863  \\
 4243 & 25896  &  4253 & 25968  &  4259 & 25994  &  4261 & 26014  &  4271 & 26091  &  4273 & 26073  &  4283 & 26133  \\
 4289 & 26178  &  4297 & 26213  &  4327 & 26421  &  4337 & 26482  &  4339 & 26491  &  4349 & 26560  &  4357 & 26651  \\
 4363 & 26656  &  4373 & 26704  &  4391 & 26864  &  4397 & 26882  &  4409 & 26967  &  4421 & 27021  &  4423 & 27043  \\
 4441 & 27173  &  4447 & 27221  &  4451 & 27230  &  4457 & 27270  &  4463 & 27306  &  4481 & 27417  &  4483 & 27422  \\
 4493 & 27512  &  4507 & 27619  &  4513 & 27641  &  4517 & 27683  &  4519 & 27678  &  4523 & 27658  &  4547 & 27877  \\
 4549 & 27845  &  4561 & 27949  &  4567 & 27958  &  4583 & 28080  &  4591 & 28144  &  4597 & 28195  &  4603 & 28173  \\
 4621 & 28322  &  4637 & 28446  &  4639 & 28421  &  4643 & 28454  &  4649 & 28487  &  4651 & 28525  &  4657 & 28527  \\
 4663 & 28614  &  4673 & 28667  &  5003 & 30823  &  6007 & 37344  &  7001 & 43831  &  8009 & 50515  &  &  \\
\hline
\end{tabular}
\end{center}

\newpage
\begin{center}
\textbf{Table 2.} Sizes $t_{2}^{G}(3,q)=t_{2}^{G}$ of
complete caps\footnote{The sizes improving the ones from \cite[Table 7]{DFMP2009} and \cite[Section 3]{Platoni}
 are written  in bold font} in $\mathrm{PG}(3,q)$ obtained using
randomized greedy algorithms, $q\in G_{3}$ \medskip\\
{\small
\begin{tabular}{@{}r@{\,\,\,}r@{\,\,}|@{\,\,}r@{\,\,\,}r@{\,\,}|@{\,\,}r@{\,\,\,}r@{\,\,}|@{\,\,}
r@{\,\,\,}r@{\,\,}|@{\,\,}r@{\,\,\,}r@{\,\,}|@{\,\,}r@{\,\,\,}r@{\,\,}|@{\,\,}r@{\,\,\,}r@{}}
\hline
$q$&$t_2^G$&$q$&$t_2^G$&$q$&$t_2^G$&$q$&$t_2^G$&$q$&$t_2^G$&$q$&$t_2^G$&$q$
&$t_{2_{\vphantom{H_H}}}^{G^{\vphantom{H^H}}}$\\
\hline
 2 & 5  &  3 & 8  &  5 & 12  &  7 & 17  &  11 & 30  &  13 & 36  &  17 & 51  \\
 19 & 58  &  23 & 72  &  29 & 96  &  31 & 104  &  37 & 128  &  41 & 145  &  43 & 153  \\
 47 & 169  &  53 & 195  &  59 & 220  &  \textbf{61} & \textbf{229}  & \textbf{67} & \textbf{255}  &  \textbf{71} & \textbf{273}  &  \textbf{73} & \textbf{282}  \\
 \textbf{79} & \textbf{309}  &  \textbf{83} & \textbf{327}  &  89 & 355  & \textbf{97} & \textbf{392}  & \textbf{101 }& \textbf{ 412  }& \textbf{103 }& \textbf{422}& \textbf{107 }& \textbf{ 437}\\\textbf{
 109 }& \textbf{447}& \textbf{113 }& \textbf{ 466  }& \textbf{127 }& \textbf{536}& \textbf{131 }& \textbf{ 560  }& \textbf{137 }& \textbf{ 583  }& \textbf{139 }& \textbf{ 598  }& \textbf{149 }& \textbf{ 646  }\\\textbf{
 151 }& \textbf{ 653  }& \textbf{157 }& \textbf{ 687  }& \textbf{163 }& \textbf{ 717  }& \textbf{167 }& \textbf{ 734  }& \textbf{173 }& \textbf{ 761  }& \textbf{179 }& \textbf{ 789  }& \textbf{181 }& \textbf{ 807  }\\\textbf{
 191 }& \textbf{ 843  }& \textbf{193 }& \textbf{ 859  }& \textbf{197 }& \textbf{ 873  }& \textbf{199 }& \textbf{ 884  }& \textbf{211 }& \textbf{ 944  }& \textbf{223 }& \textbf{ 1007  }& \textbf{227 }& \textbf{ 1028  }\\\textbf{
 229 }& \textbf{ 1043  }& \textbf{233 }& \textbf{ 1056  }& \textbf{239 }& \textbf{ 1091  }& \textbf{241 }& \textbf{ 1095  }& \textbf{251 }& \textbf{ 1155  }& \textbf{257 }& \textbf{ 1180  }& \textbf{263 }& \textbf{ 1215  }\\\textbf{
 269 }& \textbf{ 1245  }& \textbf{271 }& \textbf{ 1249  }& \textbf{277 }& \textbf{ 1286  }& \textbf{281 }& \textbf{ 1304  }& \textbf{283 }& \textbf{ 1313  }& \textbf{293 }& \textbf{ 1363  }& \textbf{307 }& \textbf{ 1441  }\\\textbf{
 311 }& \textbf{ 1459  }& \textbf{313 }& \textbf{ 1470  }& \textbf{317 }& \textbf{ 1489  }& \textbf{331 }& \textbf{ 1560  }& \textbf{337 }& \textbf{ 1597  }& \textbf{347 }& \textbf{ 1644  }& \textbf{349 }& \textbf{ 1655}  \\
 353 & 1670  &  359 & 1709  &  367 & 1748  &  373 & 1787  &  379 & 1814  &  383 & 1836  &  389 & 1866  \\
 397 & 1910  &  401 & 1929  &  409 & 1973  &  419 & 2023  &  421 & 2037  &  431 & 2090  &  433 & 2110  \\
 439 & 2134  &  443 & 2149  &  449 & 2191  &  457 & 2228  &  461 & 2252  &  463 & 2276  &  467 & 2282  \\
 479 & 2351  &  487 & 2391  &  491 & 2414  &  499 & 2463  &  503 & 2478  &  509 & 2509  &  521 & 2586  \\
 523 & 2591  &  541 & 2696  &  547 & 2720  &  557 & 2777  &  563 & 2813  &  569 & 2846  &  571 & 2857  \\
 577 & 2890  &  587 & 2949  &  593 & 2982  &  599 & 3009  &  601 & 3018  &  607 & 3042  &  613 & 3088  \\
 617 & 3114  &  619 & 3103  &  631 & 3174  &  641 & 3227  &  643 & 3243  &  647 & 3265  &  653 & 3298  \\
 659 & 3322  &  661 & 3325  &  673 & 3405  &  677 & 3419  &  683 & 3462  &  691 & 3516  &  701 & 3563  \\
 709 & 3601  &  719 & 3638  &  727 & 3702  &  733 & 3750  &  739 & 3760  &  743 & 3794  &  751 & 3841  \\
 757 & 3868  &  761 & 3893  &  769 & 3925  &  773 & 3949  &  787 & 4030  &  797 & 4095  &  809 & 4166  \\
 811 & 4181  &  821 & 4220  &  823 & 4217  &  827 & 4248  &  829 & 4264  &  839 & 4323  &  853 & 4398  \\
 857 & 4431  &  859 & 4421  &  863 & 4454  &  877 & 4529  &  881 & 4555  &  883 & 4548  &  887 & 4593  \\
 907 & 4701  &  911 & 4739  &  919 & 4768  &  929 & 4829  &  937 & 4872  &  941 & 4881  &  947 & 4921  \\
 953 & 4971  &  967 & 5037  &  971 & 5060  &  977 & 5096  &  983 & 5125  &  991 & 5169  &  997 & 5226  \\
 1009 & 5283  &  1013 & 5308  &  1019 & 5345  &  1021 & 5347  &  1031 & 5402  &  1033 & 5404  &  1039 & 5445  \\
 1049 & 5503  &  1051 & 5514  &  1061 & 5564  &  1063 & 5578  &  1069 & 5616  &  1087 & 5705  &  1091 & 5742  \\
 1093 & 5745  &  1097 & 5766  &  1103 & 5794  &  1109 & 5835  &  1117 & 5896  &  1123 & 5912  &  1129 & 5936  \\
 1151 & 6088  &  1153 & 6093  &  1163 & 6137  &  1171 & 6192  &  1181 & 6269  &  1187 & 6299  &  1193 & 6323  \\
 1201 & 6381  &  1213 & 6428  &  1217 & 6455  &  1223 & 6482  &  1229 & 6509  &  1231 & 6533  &  1237 & 6588  \\
 1249 & 6645  &  1259 & 6701  &  1277 & 6821  &  1279 & 6836  &  1283 & 6867  &  1289 & 6889  &  1291 & 6889  \\
 1297 & 6923  &  1301 & 6941  &  1303 & 6969  &  1307 & 6990  &  1319 & 7062  &  1321 & 7081  &  1327 & 7098  \\
 1361 & 7283  &  1367 & 7321  &  1373 & 7365  &  1381 & 7393  &  1399 & 7494  &  1409 & 7567  &  1423 & 7636  \\
 1427 & 7640  &  1429 & 7688  &  1433 & 7716  &  1439 & 7748  &  1447 & 7782  &  1451 & 7785  &  1453 & 7819  \\
 1459 & 7850  &  1471 & 7910  &  1481 & 7970  &  1483 & 7966  &  1487 & 8027  &  1489 & 8036  &  1493 & 8046  \\
 1499 & 8093  &  1511 & 8137  &  1523 & 8225  &  1531 & 8269  &  1543 & 8345  &  1549 & 8388  &  1553 & 8413  \\
 1559 & 8437  &  1567 & 8490  &  1571 & 8507  &  1579 & 8548  &  1583 & 8579  &  1597 & 8644  &  1601 & 8690  \\
 \hline
\end{tabular}
}
\end{center}

\newpage
\begin{center}
\textbf{Table 2. Continue.} Sizes $t_{2}^{G}(3,q)=t_{2}^{G}$ of
complete caps in $\mathrm{PG}(3,q)$ obtained using
randomized greedy algorithms, $q\in G_{3}$ \medskip\\
{\small
\begin{tabular}{@{}r@{\,\,\,}r@{\,\,}|@{\,\,}r@{\,\,\,}r@{\,\,}|@{\,\,}r@{\,\,\,}r@{\,\,}|@{\,\,}
r@{\,\,\,}r@{\,\,}|@{\,\,}r@{\,\,\,}r@{\,\,}|@{\,\,}r@{\,\,\,}r@{\,\,}|@{\,\,}r@{\,\,\,}r@{}}
\hline
$q$&$t_2^G$&$q$&$t_2^G$&$q$&$t_2^G$&$q$&$t_2^G$&$q$&$t_2^G$&$q$&$t_2^G$&$q$
&$t_{2_{\vphantom{H_H}}}^{G^{\vphantom{H^H}}}$\\
 \hline
 1607 & 8715  &  1609 & 8722  &  1613 & 8748  &  1619 & 8764  &  1621 & 8788  &  1627 & 8819  &  1637 & 8904  \\
 1657 & 9011  &  1663 & 9049  &  1667 & 9065  &  1669 & 9067  &  1693 & 9215  &  1697 & 9223  &  1699 & 9258  \\
 1709 & 9314  &  1721 & 9383  &  1723 & 9384  &  1733 & 9437  &  1741 & 9488  &  1747 & 9522  &  1753 & 9578  \\
 1759 & 9616  &  1777 & 9708    &  1783 & 9739  &  1787 & 9785  &  1789 & 9771  &  1801 & 9864  &  1811 & 9917  \\
 1823 & 9969  &  1831 & 10042   &  1847 & 10117  &  1861 & 10218  &  1867 & 10250  &  1871 & 10283  &  1873 & 10268  \\
 1877 & 10298  &  1879 & 10324  &  1889 & 10398  &  1901 & 10448  &  1907 & 10473  &  1913 & 10522  &  1931 & 10634  \\
 1933 & 10647  &  1949 & 10727  &  1951 & 10732  &  1973 & 10887  &  1979 & 10919  &  1987 & 10971  &  1993 & 11018  \\
 1997 & 11044  &  1999 & 11039  &  2003 & 11062  &  2011 & 11092  &  2017 & 11148  &  2027 & 11208  &  2029 & 11225  \\
 2039 & 11307  &  2053 & 11376  &  2063 & 11435  &  2069 & 11438  &  2081 & 11534  &  2083 & 11537  &  2087 & 11559  \\
 2089 & 11587  &  2099 & 11637  &  2111 & 11714  &  2113 & 11733  &  2129 & 11828  &  2131 & 11850  &  2137 & 11898  \\
 2141 & 11906  &  2143 & 11932  &  2153 & 11974  &  2161 & 12032  &  2179 & 12138  &  2203 & 12262  &  2207 & 12283  \\
 2213 & 12325  &  2221 & 12387  &  2237 & 12488  &  2239 & 12525  &  2243 & 12529  &  2251 & 12578  &  2267 & 12669  \\
 2269 & 12686  &  2273 & 12694  &  2281 & 12746  &  2287 & 12781  &  2293 & 12815  &  2297 & 12854  &  2309 & 12897  \\
 2311 & 12920  &  2333 & 13047  &  2339 & 13089  &  2341 & 13147  &  2347 & 13145  &  2351 & 13168  &  2357 & 13219  \\
 2371 & 13285  &  2377 & 13324  &  2381 & 13348  &  2383 & 13375  &  2389 & 13403  &  2393 & 13417  &  2399 & 13488  \\
 2411 & 13608  &  2417 & 13563  &  2423 & 13603  &  2437 & 13713  &  2441 & 13706  &  2447 & 13750  &  2459 & 13839  \\
 2467 & 13876  &  2473 & 13883  &  2477 & 13943  &  2503 & 14121  &  2521 & 14250  &  2531 & 14315  &  2539 & 14367  \\
 2543 & 14394  &  2549 & 14404  &  2551 & 14431  &  2557 & 14462  &  2579 & 14593  &  2591 & 14611  &  2593 & 14634  \\
 2609 & 14714  &  2617 & 14774  &  2621 & 14792  &  2633 & 14874  &  2647 & 14972  &  2657 & 15030  &  2659 & 15054  \\
 2663 & 15059  &  2671 & 15101  &  2677 & 15180  &  2683 & 15223  &  2687 & 15235  &  2689 & 15228  &  2693 & 15257  \\
 2699 & 15309  &  2707 & 15364  &  2711 & 15387  &  2713 & 15394  &  2719 & 15448  &  2729 & 15514  &  2731 & 15521  \\
 2741 & 15579  &  2749 & 15635  &  2753 & 15644  &  2767 & 15744  &  2777 & 15791  &  2789 & 15897  &  2791 & 15887  \\
 2797 & 15919  &  2801 & 15942  &  2803 & 15968  &  2819 & 16095  &  2833 & 16150  &  2837 & 16174  &  2843 & 16219  \\
 2851 & 16257  &  2857 & 16326  &  2861 & 16342  &  2879 & 16485  &  2887 & 16484  &  2897 & 16533  &  2903 & 16560  \\
 2909 & 16636  &  2917 & 16702  &  2927 & 16751  &  2939 & 16841  &  2953 & 16896  &  2957 & 16923  &  2963 & 16974  \\
 2969 & 16981  &  2971 & 17052  &  2999 & 17188  &  3001 & 17192  &  3011 & 17264  &  3019 & 17339  &  3023 & 17326  \\
 3037 & 17466  &  3041 & 17451  &  3049 & 17530  &  3061 & 17589  &  3067 & 17628  &  3079 & 17684  &  3083 & 17719  \\
 3089 & 17727  &  3109 & 17819  &  3119 & 17897  &  3121 & 17900  &  3137 & 17995  &  3163 & 18196  &  3167 & 18205  \\
 3169 & 18214  &  3181 & 18287  &  3187 & 18329  &  3191 & 18335  &  3203 & 18442  &  3209 & 18474  &  3217 & 18481  \\
 3221 & 18562  &  3229 & 18620  &  3251 & 18775  &  3253 & 18768  &  3257 & 18774  &  3259 & 18849  &  3271 & 18865  \\
 3299 & 19054  &  3301 & 19036  &  3307 & 19130  &  3313 & 19165  &  3319 & 19215  &  3323 & 19218  &  3329 & 19228  \\
 3331 & 19287  &  3343 & 19367  &  3347 & 19361  &  3359 & 19466  &  3361 & 19449  &  3371 & 19547  &  3373 & 19515  \\
 3389 & 19643  &  3391 & 19632  &  3407 & 19740  &  3413 & 19806  &  3433 & 19897  &  3449 & 19965  &  3457 & 20030  \\
 3461 & 20058  &  3463 & 20062  &  3467 & 20102  &  3469 & 20093  &  3491 & 20234  &  3499 & 20285  &  3511 & 20411  \\
 3517 & 20430  &  3527 & 20511  &  3529 & 20494  &  3533 & 20548  &  3539 & 20580  &  3541 & 20597  &  3547 & 20612  \\
 3557 & 20697  &  3559 & 20667  &  3571 & 20782  &  3581 & 20839  &  3583 & 20850  &  3593 & 20947  &  3607 & 20996  \\
 3613 & 21081  &  3617 & 21037  &  3623 & 21164  &  3631 & 21220  &  3637 & 21193  &  3643 & 21215  &  3659 & 21367  \\
 3671 & 21390  &  3673 & 21489  &  3677 & 21463  &  3691 & 21514  &  3697 & 21529  &  3701 & 21645  &  3803 & 22178  \\
3907 & 22859  &  4001 & 23401  &  4289 & 25225  &  &  &  &  &  &  &  & \\
 \hline
 \end{tabular}
 }
\end{center}

\newpage
\begin{center}
\textbf{Table 3.} Sizes $t_2^L(4,q)=t_2^L$ of complete lexicaps in
$\mathrm{PG}(4,q)$, $q\in L_{4}$ \medskip\\
 \begin{tabular}{@{}r@{\,\,\,}r@{\,\,}|@{\,\,}r@{\,\,\,}r@{\,\,}|@{\,\,}r@{\,\,\,}r@{\,\,}|@{\,\,}
r@{\,\,\,}r@{\,\,}|@{\,\,}r@{\,\,\,}r@{\,\,}|@{\,\,}r@{\,\,\,}r@{}}
 \hline
$q$ & $t_{2}^{L}$ &
$q$ & $t_{2}^{L}$ &
$q$ & $t_{2}^{L}$ &
$q$ & $t_{2}^{L}$ &
$q$ & $t_{2}^{L}$ &
$q$ & $t_{2_{\vphantom{H_H}}}^{L^{ \vphantom{H^H}}}$ \\
 \hline
 2 & 16  &  3 & 16  &  5 & 44  &  7 & 74  &  11 & 157  &  13 & 203  \\
 17 & 316  &  19 & 378  &  23 & 509  &  29 & 745  &  31 & 833  &  37 & 1095  \\
 41 & 1296  &  43 & 1396  &  47 & 1602  &  53 & 1937  &  59 & 2302  &  61 & 2433  \\
 67 & 2831  &  71 & 3086  &  73 & 3228  &  79 & 3681  &  83 & 3960  &  89 & 4436  \\
 97 & 5069  &  101 & 5409  &  103 & 5581  &  107 & 5920  &  109 & 6095  &  113 & 6445  \\
 127 & 7761  &  131 & 8138  &  137 & 8737  &  139 & 8943  &  149 & 9967  &  151 & 10201  \\
 157 & 10855  &  163 & 11503  &  167 & 11972  &  173 & 12620  &  179 & 13312  &  181 & 13544  \\
 191 & 14763  &  193 & 15026  &  197 & 15489  &  199 & 15755  &  211 & 17255  &  223 & 18818  \\
 227 & 19371  &  229 & 19633  &  233 & 20157  &  239 & 20985  &  241 & 21282  &  251 & 22687  \\
 257 & 23511  &  263 & 24404  &  269 & 25342  &  271 & 25588  &  277 & 26497  &  281 & 27092  \\
 283 & 27386  &  293 & 28913  &  307 & 31160  &  311 & 31754  &  313 & 32100  &  317 & 32772  \\
 331 & 35017  &  337 & 36027  &  347 & 37724  &  349 & 38090  &  353 & 38793  &  359 & 39792  \\
 367 & 41182  &  373 & 42261  &  379 & 43332  &  383 & 44038  &  389 & 45118  &  397 & 46601  \\
 401 & 47359  &  409 & 48830  &  419 & 50717  &  421 & 51132  &  431 & 52980  &  433 & 53365  \\
 439 & 54573  &  443 & 55309  &  449 & 56538  &  457 & 58157  &  461 & 58926  &  463 & 59287  \\
 467 & 60094  &  479 & 62545  &  487 & 64212  &  491 & 64999  &  499 & 66689  &  503 & 67494  \\
 509 & 68745  &  521 & 71375  &  523 & 71800  &  541 & 75708  &  547 & 77032  &  557 & 79195  \\
 563 & 80569  &  569 & 81925  &  571 & 82440  &  577 & 83791  &  587 & 86086  &  593 & 87403  \\
 599 & 88749  &  601 & 89304  &  607 & 90711  &  613 & 92127  &  617 & 93061  &  619 & 93474  \\
 631 & 96338  &  641 & 98664  &  643 & 99215  &  647 & 100114  &  653 & 101572  &  659 & 103041  \\
 661 & 103550  &  673 & 106602  &  677 & 107510  &  683 & 108942  &  691 & 111219  &  701 & 113524  \\
 709 & 115495  &  719 & 118167  &  727 & 120203  &  733 & 121660  &  739 & 123174  &  743 & 124346  \\
 751 & 126409  &  757 & 128030  &  761 & 129041  &  769 & 131161  &  773 & 132247  &  787 & 136021  \\
 797 & 138787  &  809 & 141999  &  811 & 142507  &  821 & 145352  &  823 & 145837  &  827 & 146975  \\
 829 & 147565  &  839 & 150278  &  853 & 154310  &  857 & 155338  &  859 & 156008  &  863 & 157105  \\
 877 & 161115  &  881 & 162240  &  883 & 162756  &  887 & 163911  &  907 & 169825  &  911 & 170936  \\
 919 & 173333  &  929 & 176239  &  937 & 178683  &  941 & 179900  &  947 & 181636  &  953 & 183357  \\
 967 & 187579  &  971 & 188822  &  977 & 190636  &  983 & 192545  &  991 & 194948  &  997 & 196781  \\
 1009 & 200504  &  1013 & 201779  &  1019 & 203602  &  1021 & 204253  &  1031 & 207445  &  1033 & 208062  \\
 1039 & 209945  &  1049 & 213006  &  1051 & 213705  &  1061 & 216874  &  1063 & 217563  &  1087 & 225218  \\
 1091 & 226516  &  1093 & 227235  &  1097 & 228418  &  1103 & 230322  &  1109 & 232422  &  1117 & 235049  \\
 1123 & 236956  &  1129 & 238987  &  1151 & 246325  &  1153 & 246977  &  1163 & 250167  &  1171 & 253092  \\
 1181 & 256352  &  1187 & 258447  &  1193 & 260448  &  1201 & 263240  &  1213 & 267271  &  1217 & 268644  \\
 1223 & 270732  &  1229 & 272856  &  1231 & 273487  &  1237 & 275468  &  1249 & 279794  &  1259 & 283322  \\
 1277 & 289573  &  1279 & 290388  &  1283 & 291795  &  1289 & 293821  &  1291 & 294584  &  1297 & 296679  \\
 1301 & 298264  &  1303 & 298829  &  1307 & 300381  &  1319 & 304605  &  1321 & 305314  &  1361 & 319836  \\
 1409 & 337667  &  &  &  &  &  &  &  &  &  &  \\
 \hline
 \end{tabular}
\end{center}

\newpage
\begin{center}
\textbf{Table 4.} Sizes $t_2^G(4,q)=t_2^G$ of complete
caps\footnote{The sizes improving the ones from \cite[Table 8]{DFMP2009} and \cite[Theorem 1.1]{FaPasSch_2012}
 are written  in bold font}  in $\mathrm{PG}(4,q)$, $q\in G_{4}$, obtained by
greedy algorithms\medskip \\
 \begin{tabular}{@{}r@{\,\,\,\,}r@{\,\,}|@{\,\,}r@{\,\,\,\,}r@{\,\,}|@{\,\,}r@{\,\,\,\,}r@{\,\,}|@{\,\,}
r@{\,\,\,\,}r@{\,\,}|@{\,\,}r@{\,\,\,\,}r@{\,\,}|@{\,\,}r@{\,\,\,\,}r@{\,\,}|@{\,\,}r@{\,\,\,\,}r@{}}
 \hline
$q$ & $t_{2}^{G}$ &
$q$ & $t_{2}^{G}$ &
$q$ & $t_{2}^{G}$ &
$q$ & $t_{2}^{G}$ &
$q$ & $t_{2}^{G}$ &
$q$ & $t_{2}^{G}$ &
$q$ & $t_{2_{\vphantom{H_H}}}^{G^{\vphantom{H^H}}}$ \\
 \hline
 2 & 9  &  3 & 11  &  5 & 31  &  7 & 56  &  11 & 121  &  13 & 162  &  \textbf{17} & \textbf{254}  \\
 19 & 309  &  23 & 425  &  \textbf{29} & \textbf{625}  &  31 & 695  &  37 & 935  &  41 & 1106  &  43 & 1194  \\
 47 & 1386  &  \textbf{53} & \textbf{1687}  &  59 & 2013  &  61 & 2123  &  \textbf{67} & \textbf{2476}  &  71 & 2723  &  \textbf{73} & \textbf{2870}  \\
 79 & 3253  &  \textbf{83} & \textbf{3535}  &  \textbf{89} & \textbf{3982}  &  \textbf{97} & \textbf{4526}  &  \textbf{101} & \textbf{4868}  &  \textbf{103} & \textbf{5023}  &  \textbf{107} & \textbf{5326}  \\
 \textbf{109} & \textbf{5512}  & \textbf{ 113} & \textbf{5814}  &  \textbf{127} & \textbf{7021}  &  131 & 7437  &  137 & 7987  &  139 & 8161  &  149 & 9118  \\
 151 & 9316  &  157 & 9899  &  163 & 10510  &  167 & 10958  &  173 & 11545  &  179 & 12223  &  181 & 12451  \\
 191 & 13573  &  193 & 13798  &  197 & 14266  &  199 & 14511  &  211 & 15902  &  223 & 17360  &  227 & 17876  \\
 229 & 18162  &  233 & 18605  &  239 & 19382  &  241 & 19682  &  251 & 20997  &  257 & 21766  &  263 & 22577  \\
 269 & 23453  &  271 & 23702  &  277 & 24532  &  281 & 25122  &  283 & 25391  &  293 & 26821  &  307 & 28892  \\
 311 & 29490  &  313 & 29785  &  317 & 30380  &  331 & 32529  &  337 & 33513  &  347 & 35118  &  349 & 35423  \\
 353 & 36082  &  359 & 37065  &  367 & 38371  &  373 & 39422  &  379 & 40383  &  383 & 41100  &  389 & 42123  \\
 397 & 43984  &  401 & 44547  &  409 & 45964  &  419 & 47636  &  421 & 47899  &  431 & 49819  &  433 & 50206  \\
 439 & 51320  &  443 & 52148  &  449 & 53257  &  457 & 54478  &  461 & 55632  &  463 & 56057  &  &  \\
 \hline
 \end{tabular}
\end{center}

\end{document}